# Simplicial structures in ecological networks

Udit Raj[1], Shashankaditya Upadhyay[2], Moumita Karmakar[3], Sudeepto Bhattacharya[1] *

[1] Department of Mathematics, School of Natural Sciences, Shiv Nadar University Delhi-NCR, Gautam Buddha Nagar, Greater Noida 201314 Uttar Pradesh, India

[2] Department of Physical Sciences, Indian Institute of Science Education and Research Kolkata, Mohanpur 741246 West Bengal, India.

[3] Centre for Public Affairs and Critical Theory (C-PACT), School of Humanities and Social Sciences, Shiv Nadar University Delhi-NCR, Gautam Buddha Nagar, Greater Noida 201314 Uttar Pradesh, India

∗ Corresponding author. E-mail: sudeepto.bhattacharya@snu.edu.in

**Abstract**

An ecological network is a formal representation of a specific type of interaction in a corresponding ecosystem. Such networks have traditionally been modelled as encoding exclusively pairwise interactions among the fundamental units of ecosystems and have been represented and analysed using graph-theoretic methods. However, many real-world ecosystems may entertain non-binary, polyadic relations between their units, which cannot be captured by the pairwise interaction methods, but require higher-order interaction framework, and consequently the corresponding ecological networks cannot be modelled using graph-theoretic framework. This work gives a structural definition of ecological network suitable for modelling all orders of interactions between the fundamental units of the corresponding ecological system, including and going beyond the pairwise interaction framework. Carbon mediation between units of some select ecosystems are studied by modelling the corresponding ecological networks as simplicial complexes following the definition. The concept of graph centrality measure has been extended to simplicial centrality, and some important centrality measures of these networks at various structural levels of the complexes have been calculated. The centrality measures reveal valuable structural information including information about those vertices that are more likely to participate in higher-order interactions, as well as inform whether there is a difference in the ranks of vertices for these higher-order networks based on graph centrality and simplicial centrality measures.

# 1 Introduction

An ecosystem is structurally an assemblage of the interacting ecological units that comprise it. The interactions among the units are encoded as relations between the units, and are primarily responsible for the emergence of collective behaviour in the system, qualifying every ecosystem as a complex adaptive system [1, 2, 3]. An ecosystem may typically comprise a range of relations, from the binary ones that encode pairwise interactions among the ecological units, to more general *m*-ary ones that encode non-pairwise interactions within a collection of more than two units. Each of the relations constituting a given ecosystem can be mathematically modelled and represented by a corresponding ecological network. Formal representation of ecosystems as ecological networks has evolved as a very powerful and an equally successful framework, and over the past two decades has been instrumental in facilitating important insights and understanding of the structural organisation of ecosystems, and the observed processes and properties of the ecosystem due to a specific structure [4 – 22].

A network, as formally defined in [22], can be used to quantitatively model empirical data of pairwise interaction for an ecosystem. A graph is traditionally used as the underlying combinatorial object, and a network is modelled based on this graph with the edges depicting the interactions between vertices, which represent the fundamental units of the ecosystem [23 - 27].

Majority of the ecological networks proposed and studied in the scholarship, however, are premised on the framework of simple graphs as mentioned above, and therefore are based on the assumption that the complex interactions and the structural organisation of a given ecosystem can be modelled exclusively by the pairwise interactions among the units of the system. Such a graph theory-based ecological network formalism for an ecosystem properly models the represented system provided the system admits only pairwise interactions among its fundamental units. However, for ecosystems that may admit not only pairwise interactions encoded as binary relations, but also non-pairwise interactions and thus non-binary relations among their units, such representations face a twin mathematical challenge:(i) due to the definition of graph as a combinatorial object having every edge defined as a set of pair of its vertices, the ecological network modelled with graph-theoretic framework is necessarily limited to only capture the pairwise interactions and therefore the binary relations existing in the system, and hence is successful in capturing the system interactions only partially, and therefore, (ii) there exists a possibility of masking those complex behaviours of the modelled system that are emergent due to the existing non-pairwise interactions among the units of the system [28].

Often these two constraints render the modelled ecological network as an oversimplified representation of a corresponding ecosystem, thereby failing to capture the structural intricacies resulting from non-pairwise interactions of the units of the ecosystem. The actual problem for this failure, as explained above, lies with the definition of ecological network, which is based on a graph theoretic framework that only allows for pairwise interactions to be considered. It is therefore of an overriding need that mathematical framework beyond the graph theoretic one be adopted to define ecological network so that it is structurally able to capture and represent all orders of interactions in the corresponding ecosystem [28]. The principal objective of this work is to address the first of the challenges mentioned earlier, and propose a definition for ecological network using the framework of algebraic topology as an object essentially comprising a simplicial complex, so as to account for and represent all orders of interactions in the modelled ecosystem.

As illustrations of our stated objective, the work applies the constructed definition to some select ecosystems focusing on the relation of carbon mediation, in which it has been assumed that the relation between the vertices (ecological units, biotic as well as abiotic) of the networks is due to carbon mediation [29,30]. However, as previously stated, all real-world systems may not necessarily have purely dyadic or pairwise relations among their vertices. Rather, the system's relations may be polyadic in nature, involving the interaction of more than two vertices as a collection, and this aspect of interaction cannot be overlooked, which will be referred to as higher-order interaction in this work [31,32]. This point could be appreciated using some examples of mutualism and also in the instance of a recent study, according to which the genus *Anolis* has a large number of ectomorphs. *Anolis sp.* (A) feeds on both crickets (B) and insects (C), and crickets eat small insects as well (as they prefer carnivorous diet). As a result, A, B, and C interact as a set of three units, and therefore share a non-binary relation of carbon mediation among them, exemplifying higher-order interactions in ecosystems [ 33, 34]. The present work, through defining the concept of ecological network, argues that these non-pairwise interactions can be appropriately captured by the combinatorial object called simplicial complex in algebraic topology.

The structural differences between carbon mediation networks by exclusively considering the pairwise interactions have been established previously by random projections in $\mathbb{R}^2$ and $\mathbb{R}^3$ of weighted spectral distribution (WSD) for certain ecosystems taken from different environments [35]. However, while comparing these networks on the basis of structure, the polyadic relations encoding the non-pairwise interactions of the ecological units must be kept in consideration and the construction of an ecological network representing the interactions as a simplicial complex following our definition for the same would be an efficient method to understand the interactions better [36, 37].

For a comprehensive insight into the structural differences at different structural levels of an ecological network, centrality measures such as betweenness centrality, closeness centrality, degree centrality, subgraph centrality, and eigenvector centrality must be calculated. In this work, we propose to calculate these metrics for higher-order ecological networks by generalising the definitions of these measures to non-pairwise, higher-order interactions from the graph-theoretic pairwise ones. Simplicial complexes allow us to examine the structural intricacies of these networks in detail, allowing us to gain a better understanding of carbon mediation networks. On the one hand, the structures of these networks are frequently thought to be similar [38, 39, 40], and higher-order interactions in ecosystems have been postulated or considered to be extremely rare by some researchers, some recent research articles, on the other hand, have found evidence of a large number of higher-order interactions in ecosystems, which can help us better comprehend the structural complexity and the dynamics of ecosystem processes [41, 42]. As mentioned above, we illustrate the applicability of our constructed definition of ecological network by considering the carbon mediation networks in some select ecosystems. We study and structurally analyse these networks as a simplicial complexes, and compare the graph and simplicial centralities for the constituent simplices of dimension 0, 1 and 2 in these networks, with the measures giving us the comparisons of structures at the respective simplicial dimensions.

In this study while focusing on the non-pairwise interactions, we pose some questions: first, is it possible to construct simplicial complex of each of the network studied in this work? Second, is simplicial complex of a connected graph $G$ always connected? Third, is there any difference in the rankings of vertices on the basis of graph centrality and simplicial centrality measures for carbon mediation networks at different structural levels? The present work, along with addressing the stated objectives, also addresses these important structural questions.

The carbon mediation networks explored in this work have been designated as aquatic interaction networks if they originate from aquatic ecosystems, while the ones designated as terrestrial networks are embedded in terrestrial ecosystems. In order to explore the interactions in these networks and thus reach our stated objective, we begin by constructing a simplicial complex from each of the graph corresponding to a given network. By defining the corresponding adjacency matrices for each complex, the simplicial centralities are determined at the edge (1-simplex) and triangle (2-simplex) levels, taking into account the interactions among the 1-simplices and 2-simplices of the higher-order networks. The specific centrality indices are then calculated using the simplicial centrality scores as markers to see if there is a difference in the vertices ranking of carbon mediation networks emanating from the two groups of networks. Cliques of various sizes have been calculated to assist us in creating the clique complex of each network to acquire a better understanding of the structure of the networks [43, 44].

The next section provides the necessary mathematical preliminaries required to advance our arguments. Definitions 1 to 11 presented in this section are found in the literature [45, 46, 47, 48, 49, 50, 51] and have been repeated here to make the work self-contained. We begin with the definition of ecological systems followed by that of the simplicial complex, and end the section by defining ecological network in Definition 12. In section 3, we employ this definition to construct the corresponding ecological networks obtained from publicly available data, while again repeating some of the required definitions from the literature. Section 4 reports the results obtained from our computations, while Section 5 concludes the work with a short discussion on the modelling.

## 2 Preliminaries

Some of the definitions presented in this section are in terms of structure, and are essentially combinatorial. Detailed discussion on structure as an object could be found in [27, 45, 46, 47].

***Definition* 1:** An *ecological system* $\Sigma$ is an object given by the structure $\Sigma = (V, E, C)$, where $V$, the *vertex set* is a non-empty set, the universe of the system, whose elements are the combinatorial objects called *vertices* and represent the ecological units of $\Sigma$. $E$ is the set of relations, whose elements encode the relationships among the vertices, and are called the *edges*.

The elements of $E$ are denoted by $e$ and are the members of the power set of $V$, $\mathcal{P}(V)$, the set of all subsets of $V$ including the empty set $\emptyset$ and $V$ itself. That is, the elements of $E$ are subsets of $V$ of arbitrary cardinalities over $\mathbb{N}$, and $E$ thus comprises the set of relations of arbitrary arity on $V$. The set $C$ comprises the constants of the system, and may be considered as the set of all 0-ary or nullary maps on $\Sigma$. If there is no possibility of confusion and no particular need, we may not explicitly mention the set $C$ henceforth.

Definition 1 is a generalized definition of ecological system which allows modelling all relations of heterogeneous artiry that may exist in the system. For our entire discussion in what follows, we shall assume the sets $V$ and $E$ to be finite, and write $|V| = m$, $m \in \mathbb{N}$, etc, where $|V|$ denotes the number of elements in the set $V$.

***Definition 2***: Let $\Sigma = (V, E, C)$ be an ecological system. Each element of $E$ is called an *interaction* in $\Sigma$. Thus, an interaction $e = \{v_0, \dots, v_{k-1}\}$ containing $k$ number of vertices, $k \in \mathbb{N}$, $k \leq |V|$, is an element of the power set of $V$: $e \in \mathcal{P}(V)$. Therefore, $\forall e \in E, e \in \mathcal{P}(V)$, implying that $E \subseteq \mathcal{P}(V)$. For the rest of this work, we shall assume every interaction to be non-empty.

***Definition 3***: Given the above system, the *order* of an interaction $e \in E$ with $|e| = k$ is defined to be $k - 1$, and $e$ is called a $k - 1$ interaction.

Interactions are classified as *higher-order* if $|e| = k \geq 3$, and as *lower-order* if $|e| = k \leq 2$ [48]. Since each element of $E$ is an interaction in the system, therefore the $n$- ary relations that comprise $E$ encode all interactions in $\Sigma$. Owing to this fact, we shall use the term relation interchangeably with the term interaction that this relation encodes, in rest of this article. Further, if a system $\Sigma$ is such that given any interaction $e \in E$, all non-empty subsets $e'$ of $e$ also belong to $E$, then the system is said to exhibit subset dependency [49].

**2.1 Simplicial Complex**

In this sub-section, a combinatorial object known as simplicial complex in classical algebraic topology is defined and described and an argument is presented that a simplicial complex is equipped with an appropriate structure to capture interactions of all possible (finite) orders between the fundamental units of a given system. Rich literature text has been devoted to the study of simplicial complexes and the definitions that follow this paragraph to make the paper self-sufficient are essentially found in the standard texts [50,51,52]. In this work we are interested to combinatorially model carbon mediation networks as simplicial complexes, and hence shall remain primarily concerned with the concept of an abstract simplicial complex defined in the following, while making a cursory reference to its geometric counterpart to facilitate visualisation of the networks obtained.

***Definition 4***: Let $V$ be a non-empty, finite set with $|V| = n + 1, n \in \mathbb{N}$, whose elements are the vertices, and are denoted by $v_i$, $i = 0, \ldots, n$. Any member of the power set of $V$, $\mathcal{P}(V)$, is a combinatorial object called a *simplex* over $V$ (or, a simplex).

For the rest of our discussion, we shall assume that all simplices are non-empty. A simplex therefore, is a non-empty, finite subset of $V$.

The dimension of a simplex $\sigma_i = \{v_0, \ldots, v_k\}$ with $|\sigma| = k + 1$ is defined to be $|\sigma| - 1$, that is $\dim(\sigma_i) = k$. Often, $\sigma_i$ is then called a $k$–simplex. Any subset of a simplex $\sigma_i$ is called a *face* of $\sigma_i$. If $\sigma_j$ is a face of $\sigma_i$, then it is written as $\sigma_j \leqslant \sigma_i$. A 0-face (0-dimensional subset of a simplex) is called a vertex and written as $\{v_i\}$, and a 1-face is called an edge, written as $\{v_i, v_j\}$, with the indices belonging to some index set [47].

***Definition 5***: Let $V$ be the set of vertices as given in the above definition, a family $\Delta$ of simplices is called a *simplicial complex (*or *abstract simplicial complex)* if it is closed under inclusion, that is, under taking of (finite, non-empty) subsets (often also referred to as the hereditary property or the property of downward-closedness).

Though the collection $\Delta$ may either be finite or infinite, we shall assume it to be finite for our discussion.

The above definition means that the following condition holds for $\Delta$:

$\forall \sigma_i \in \Delta, (\sigma_j \subset \sigma_i \Rightarrow \sigma_j \in \Delta)$. That is, every face of a simplex in $\Delta$ must also belong to $\Delta$.
(1)

This condition further implies that two arbitrary simplices in the simplicial complex $\Delta$ are either disjoint, or they intersect as a face of $\Delta$. That is, $\forall \sigma_i, \sigma_j \in \Delta \Rightarrow (i)\; \sigma_i \cap \sigma_j = \emptyset$, or $(ii)\; \sigma_i \cap \sigma_j \in \Delta$.

The *dimension* of a simplicial complex $\Delta$, denoted by $dim(\Delta)$, is defined to be $r \geq 0$ where $r$ is the largest natural number such that $\Delta$ contains an $r$-simplex. If $dim(\Delta) = d$, then every face of dimension $d$ is called a *cell*, while that of dimension $d$-1 is called a *facet*.

***Definition 6***: Let $\Delta$ be a simplicial complex. The *boundary* of $\Delta$, denoted by $\partial\Delta$, is the set of all the faces of $\Delta$. $\partial\Delta := \{\sigma_j \mid \sigma_j \leqslant \sigma_i \in K\}$, where $\sigma_i$ is a unique $k$-simplex of $\Delta$.

The above definition implies that for a given $k$-simplex in $\Delta$, every face that is a $(k$-1$)$-simplex is a boundary face or a facet. Further,

***Definition 7***: Let $\Delta$ be a simplicial complex on a set V. Let $p \in \mathbb{N}$. The simplicial complex given by $\Omega := \{\sigma_i \in \Delta \mid dim(\sigma_i) \leq p\}$ is defined to be the p-skeleton of $\Delta$.

$\Omega$ is therefore the collection of all faces of $\Delta$, that have a dimension at most $p$.

For example, a simple graph $G = (V, E)$ with a finite, non-empty vertex set $V$ and an edge set $E$ is a 1-skeleton comprising all 0-simplices (vertices) and 1-simplices (links) of a simplicial complex $\Delta$ with $dim(\Delta) \geq 1$. Therefore, all simple graphs are simplicial complexes of dimension at most 1.

As has been defined earlier, a simplicial complex is a combinatorial object: it is a collection of simplices satisfying condition (1). Every simplicial complex $\Delta$, however, corresponds to a geometric object which is a subspace of the *m*-dimensional Euclidean space $\mathbb{R}^m$ via a mapping $\varphi: V \rightarrow \mathbb{R}^m$, where $V$ is the vertex set on which $\Delta$ is defined. The *geometric realisation* of $\Delta$ with respect to the mapping $\varphi$ is defined to be the set $|\Delta|_\varphi = \cup_{\sigma \in \Delta} |\sigma|_\varphi$. The object on the right-hand side of the equality sign is the *geometric simplex* and is a subspace of $\mathbb{R}^m$. This subspace has the set of $\varphi$ images of the vertices of the simplex $\sigma \in \Delta$ as a basis, and is thus the convex hull of points in the vertex set of $|\sigma|_\varphi$. The left-hand side object therefore becomes a subspace of $\mathbb{R}^m$ and is a topological space by virtue of the definition, called the polyhedron. Thus, the topological object $|\Delta|_\varphi$ often written as $|\Delta|$ when the mapping is understood, is the geometric (and topological) counterpart of the combinatorial object $\Delta$. Note that if the mapping $\varphi$ is an affine embedding, then the topology of $|\Delta|_\varphi$ is independent of $\varphi$. However, as our interest is in modelling the system through its combinatorial information, we shall largely deal with the abstract simplicial complex as defined in Definition 5, we shall not lay much emphasis on the geometrical counterpart. For details on the geometric realisation, the reader may refer to [50, 51]. It may be noted that simplicial complexes are mathematically equipped to interface between combinatorial (discrete) and geometric (continuous) objects, and this may be considered as one of the prime reasons for modelling higher-order interactions as simplicial complexes that would offer an understanding and visualisation of the geometric configurations that they represent [53].

***Definition 8:*** let $A = \{a_0, a_1, \ldots, a_k\}$ be a geometrically independent set of points in $\mathbb{R}^n$, $n \geq k$. Then the *k*-dimensional geometric simplex or *k*-simplex spanned by the set $A$, is the set of all those points $x \in \mathbb{R}^n$ such that

$x = \sum_{i=0}^{k} \alpha_i a_i$ where $\sum_{i=0}^{k} \alpha_i = 1$ and $\alpha_i \geq 0$

For each $i = 0,1,2,\ldots,k$   [51]

By the Definition 4, it is clear that a *k*-simplex contains $k + 1$ vertices. So, every clique of size $k$ in a graph is a (k-1)-simplex, for example a line is a (2-1) simplex, a triangle is a (3-1) simplex etc. The below table gives the geometric interpretation of simplex in $\mathbb{R}^n$.

| Simplex | Geometric interpretation in $\mathbb{R}^n$ |
|---|---|
| 0-simplex: $\langle v_i \rangle$, a single vertex | $\{v_i\}$, a single point |
| 1-simplex: $\langle v_i, v_j \rangle$, two vertices | $\{v_i, v_j\}$, an unordered pair of points determining a closed line segment joining them |
| 2-simplex: $\langle v_i, v_j, v_k \rangle$, three vertices | $\{v_i, v_j, v_k\}$, an unordered triple of points determining a triangle together with its interior, joining the three points |
| ⋮ | |
| k-simplex: $\langle v_0, \ldots, v_k \rangle$, $k+1$ vertices | $\{v_0, \ldots, v_k\}$, an unordered collection of $(k+1)$ points that determine a $(k+1)$-gon, including its interior. |

## 2.2 Adjacency in simplicial complex

Adjacency in a carbon mediation network is defined by the carbon flow between ecological units; if there is mediation of carbon between two ecological units, they are considered adjacent; otherwise, they are considered not adjacent. However, in comparison to graph theory, adjacency in the simplicial complex is a little more difficult to define. A simplicial complex has simplices of various dimensions, making it difficult to create a generalised object such as a graph adjacency matrix. Adjacency matrix in a simplicial complex $\Delta$ can be defined at different levels $k = 0, 1, 2, \ldots$ corresponding to the dimension of the simplex in $\Delta$. $k = 0$ represents the adjacency in graphs and k=1 represents the adjacency of 1-simplices.

Adjacency in simplicial complex at each level can be defined in two ways (i) lower adjacency matrix and (ii) upper adjacency matrix. $A_{low}^k$ denotes the lower adjacency matrix at $k^{th}$ level and $A_{up}^k$ denotes the upper adjacency matrix at $k^{th}$ level.

***Definition 9:*** Let $\sigma_i$ and $\sigma_j$ be two $k$-simplices. Then, the two $k$-simplices are lower adjacent if they share a common face of dimension (k-1). That is, for two distinct $k$-simplices $\sigma_i = \{w_0, w_1, \ldots, w_k\}$ and $\sigma_j = \{v_0, v_1, \ldots, v_k\}$ then $\sigma_i$ and $\sigma_j$ are lower adjacent if and only if there is a $(k\text{-}1)$-simplex $\sigma_k = \{x_0, x_1, \ldots, x_{k-1}\}$ such that $\sigma_k \subset \sigma_i$ and $\sigma_k \subset \sigma_j$. We denote lower adjacency by $\sigma_i \smile \sigma_j$.

For k=0, lower adjacency of two 0-simplices is not defined by this definition, because for two 0-simplices to be lower adjacent they should share a common face of dimension -1, which is not possible.

$$[A_{low}^k]_{ij} = \begin{cases} 1 & \text{if } \sigma_i \smile \sigma_j \\ 0 & \text{otherwise or } i = j \end{cases}$$

***Definition 10***: Let $\sigma_i$ and $\sigma_j$ be two $k$-simplices. Then, the two $k$-simplices are upper adjacent if they both are faces of the same $(k+1)$-simplex. That is, for $\sigma_i = \{w_0, w_1, \ldots, w_k\}$ and $\sigma_j = \{v_0, v_1, \ldots, v_k\}$ then $\sigma_i$ and $\sigma_j$ are upper adjacent if and only if there is a $(k+1)-$simplex $\sigma_k = \{x_0, x_1, \ldots, x_{k+1}\}$ such that $\sigma_i \subset \sigma_k$ and $\sigma_j \subset \sigma_k$. We denote the upper adjacency by $\sigma_i \frown \sigma_j$

For k=0, two 0-simplices are called to be adjacent if they are upper adjacent only, which describe the adjacency in graphs.

$$[A_{up}^k]_{ij} = \begin{cases} 1 & \text{if } \sigma_i \frown \sigma_j \\ 0 & \text{otherwise or } i = j \end{cases}$$

The following considerations arise due to the above definitions:
1. If we use lower adjacency definition as standard definition for defining adjacency in simplicial complex, then adjacency at $k=0$ level cannot be defined. Because two 0-simplices must share a negative dimension simplex in order to be lower adjacent, and negative dimension simplices are not defined.
2. If we use upper adjacency definition as standard definition for defining adjacency in simplicial complex, then adjacency between two highest dimensional simplices cannot be defined. For $k=n$ which is the highest level of any simplicial complex, the upper adjacency of two n-simplices cannot be defined.
3. If a simplicial complex has several 3-simplices but no 4-simplices, the upper adjacency of 3-simplices cannot be determined, and these 3-simplices are not adjacent to each other, so the adjacency matrix constructed will not give the required information about the network.

Each of the above points produces structural errors in constructing the network, thereby causing a loss in combinatorial information about the structure of the constructed network. In order to mitigate this loss, we adopt the following definition of the adjacency matrix for our purpose [53, 54]

***Definition 11***: Let $\sigma_i$ and $\sigma_j$ be two $k$-simplices in a simplicial complex $\Delta$. Then, for $k \geq 1$ the adjacency matrix $A^k$ at the k-level in the simplicial complex are considered adjacent if they are both lower adjacent and not upper adjacent.
For $k=0$, two 0-simplices are called to be adjacent if they are upper adjacent only, which describe the adjacency in graphs.

$$[A^k]_{ij} = \begin{cases} 1 & \text{if } \sigma_i \smile \sigma_j \text{ and } \sigma_i \not\frown \sigma_j \\ 0 & \text{if } i = j \text{ or } \sigma_i \not\smile \sigma_j \text{ or } \sigma_i \frown \sigma_j \end{cases}$$

This definition eliminates the majority of the concerns raised above, which cause us to lose vital structural information about any simplicial complex, and it also allows us to define the adjacency of higher- and lower-order simplices. Further, it enables us to calculate the relationships between a simplex's centralities and its faces, which we are specifically interested in at the vertex, edge and triangle levels.

***Definition 12***: An *ecological network* is a combinatorial object whose structure is given by $N_E = (\Delta, \Lambda)$ along with an algorithm $A$ such that for $\Lambda \neq \emptyset$, $\alpha \in \Lambda \subset \mathbb{N}$, $\Delta$ is a simplicial complex over $V$, the universe of an ecological system, with $\dim(\Delta) \geq 1$ given by the algorithm $A(\alpha)$. A network is called static if the temporal component $\Lambda$ is singleton, otherwise the network is a dynamic network.

The above definition implies that in the instance of an ecological system admitting interactions at most of order1 (that is, only lower-order, at most binary interactions), the corresponding network $N_E$ at every step indexed by $\alpha$ can be identified with a 1-skeleton of $\Delta$. Such a network is often called the underlying graph of $\Delta$, and is given by the structure $G = (V, E)$, consisting of the vertex set $V$ a non-empty set of finite abstract combinatorial objects called vertices, and a family $E$ of two-element subsets of $V$: $E \subseteq \mathcal{P}(V)$ such that $\forall e \in E, |e| = 2$, that is, each $e$ is a 1-simplex called an edge of $G$, as mentioned earlier. It is this mathematical object that gets referred to as ecological networks conventionally. Due to the overarching nature of its definition, we shall refer to the defined ecological network as higher-order network, while it would represent even the lower-order interactions that may be present in the system.

## 2.3 Simplicial centralities

The centrality indices of any network are the most effective tools for discovering its structural aspects [55]. A centrality index is a numerical value that represents the importance of a vertex in terms of its structural location in the network [56, 57]. The current work tries to apply this concept to the simplicial complex in order to better understand the relevance of a simplex of a particular order in a simplicial complex. With the help of the adjacency matrices $A^k$ corresponding to every simplex level $k$ in a simplicial complex, important $k$-simplices can be identified. For $k=0$, centrality indices give us the dominating 0-simplices of the network while for $k \geq 2$, it gives those $k$-simplices which dominantly participated in higher order interactions.

A structural index is often calculated to quantify structural characteristics or properties of a network, and every centrality index discussed in this paper is a structural index. As a result, we will start with a definition of structural index.

Before defining the structural index(s) we will see the concept of isomorphism in simplicial complex [58].

Two simplicial complexes $\Delta_1$ and $\Delta_2$ said to be isomorphic to each other if there exist simplicial maps f and g such that
$f: \Delta_1 \rightarrow \Delta_2$ and $g: \Delta_1 \rightarrow \Delta_2$ , $I_{\Delta_1} = g \circ f$ , $f \circ g = I_{\Delta_2}$ ,
Where $I_{\Delta_1}$ and $I_{\Delta_2}$ are the identity maps of $\Delta_1$ and $\Delta_2$ respectively.
Such that $f$ and $g$ are inverses of each other and are called invertible. Which implies $f$ and $g$ both are one-one and onto maps.

***Definition 13***: Let $\Delta_1$ and $\Delta_2$ are two simplicial complexes defined on a vertex set V and φ be the isomorphism between $\Delta_1$ and $\Delta_2$. A real valued function $S$ is called structural index if and only if the following condition is satisfied: $\forall \sigma \in \Delta_1: \Delta_1 \cong \Delta_2 \Longrightarrow S_{\Delta_1}(\sigma) = S_{\Delta_2}(\varphi(\sigma))$, where $S_{\Delta_1}(\sigma)$ is the value of $S(\sigma)$ in $\Delta_1$.

***Theorem 1:*** Simplicial centrality indices are invariant under isomorphic maps.
***Proof:*** Let $\Delta_1 \cong \Delta_2$ and $\varphi$ is the isomorphism map between $\Delta_1$ and $\Delta_2$.
Let $C$ be any centrality index for $\Delta_1$ such that $C: \Delta_1 \to \mathbb{R}$. A centrality index $C$ gives scores to every $k$-simplices of $\Delta_1$ to make them comparable. We need to show that same ranking of $k$-simplices is also present in $\Delta_2$.

Let $\sigma \in \Delta_1$ is a k-simplex, such that image of $\sigma$ under $C$ is $(\sigma) = x$, where $x$ is a real number, but the image of $\sigma$ under $\varphi$ is $\varphi(\sigma)$, which implies $C(\varphi(\sigma)) = x$, but $\varphi(\sigma) \in \Delta_2$. Hence there exist a $k$-simplex $\varphi(\sigma)$ in $\Delta_2$ which have the same centrality score as $\sigma$. This shows that the ranking of k-simplices under $C$ is preserved in the isomorphism map. ∎

The definition of a structural index expresses the natural requirement that a centrality measure must be invariant under isomorphism. A centrality measure needs to be a structural index, by which it can induce at least semi-order on the simplices set $\Delta$ by mapping them on an ordered set [55].

***Definition 14***: A real valued function $C$ from a simplicial complex $\Delta$ to $\mathbb{R}$, is said to be simplicial centrality index if it follows the given conditions Let $C: \Delta \to \mathbb{R}$ and $\forall \sigma_i, \sigma_j, \sigma_k \in \Delta$,

Reflexivity   : $C(\sigma_i) \leq C(\sigma_i), \forall \sigma i \in \Delta$
Comparability: $C(\sigma_i) \leq C(\sigma_j)$ or $C(\sigma_j) \leq C(\sigma_i)$
Transitivity  : $C(\sigma_i) \leq C(\sigma_j)$ and $C(\sigma_j) \leq C(\sigma_k) \Longrightarrow C(\sigma_i) \leq C(\sigma_k)$
Antisymmetry: $C(\sigma_i) \leq C(\sigma_j)$ and $C(\sigma_j) \leq C(\sigma_i) \Longrightarrow C(\sigma_i) = C(\sigma_j)$

Let $\sigma_i, \sigma_j \in \Delta$, If $C(\sigma_i) = a$ and $C(\sigma_j) = b$, as we know that $a$ and $b$ are real numbers and the set of real numbers $\mathbb{R}$ is a totally ordered set with respect to total order relation, denoted by $\leq$. If $b \leq a$, then we can say that $\sigma_i \in \Delta$ is at least as central as $\sigma_j \in \Delta$ with respect to a given centrality $C$ such that $C(\sigma_j) \leq C(\sigma_i)$.

## 3 Materials and methods

### 3.1 Construction of simplicial Complex

In this study ten ecosystems have been studied by randomly choosing five each from the terrestrial environment aquatic environment [58 - 60]. The vertices of the ecosystems are the ecological units, while they interact to mediate carbon within the system, defining the edges. In the corresponding ecological network, if two vertices in a carbon mediation network have some carbon mediation between them, they are termed adjacent vertices. Similarly, at the 1-simplex and 2-simplex level, the adjacency relation in simplicial complexes can be established. Adjacency in the simplicial complex is defined in Definition

11, which states that two 1-simplices(edges) are adjacent to each other if they have one common 0-simplex(vertex) but both 1-simplices(edges) must not participate in carbon mediation to form a 2-simplex (triangle).

For this study, the data is obtained from the Globalwebdb database website at the University of Canberra [38], which contains more than 360 food web matrices (Excel data files) and more than 120 reference papers (PDF files). It provides data in the form of adjacency matrices, and adjacency matrices for food web networks in both the aquatic and terrestrial environments have been downloaded [61, 62, 63]. Due to the small size of these adjacency matrices, the reconstruction of this data for the generation of the simplicial complex was done manually, without employing an algorithm while we remain aware of the non-trivialness and the challenges associated with the reconstruction of a general simplex from graph data. Following Definition 12, these networks were formulated as simplicial complexes, and various centrality measures at the 0, 1, and 2-simplex levels of the complex were computed.

Clique complex, Vietoris-Rips complex, and Cech complex are the three main methods for building a simplicial complex from a given network. A metric space $(X, d)$ was necessary to determine the distance between the vertices of any graph in the Vietoris-Rips and Cech complexes. Since our focus in this work is essentially combinatorial without assuming any background geometric space, we have used the clique complex technique for this work.

Further, in order to generate a simplicial complex from a given network, the network must be simple, which means the graph must not have multiple edges or loops. This can be thought of as a prerequisite for creating the simplicial complex that corresponds to the given graph. For example, the Lake-1 network, which has 20 vertices (0-simplices) and 3 loops is not simple. And simplicial adjacency matrix calculation from Definition 11 for non-simple graphs is challenging. Hence, this network first must be transformed into a simple network by deleting the loops in order to generate the clique complex. We use two graph theoretic operations to achieve this: edge deletion and vertex deletion.

On the one hand an induced subgraph of a graph can be obtained via vertex deletion, by selecting a subset of vertices that are not involved in loops, and hence resulting into a simple subgraph. However, we may lose vital information in this process because it deletes all the edges that are adjacent to the vertex that we have deleted. On the other hand, the operation of edge deletion would only delete those edges whose appears as loops. The resultant subgraph will be a simple graph and with expectedly with less information lost. Thus, the method of edge deletion was chosen to obtain simple subgraphs for this work.

**Clique Complex**

A clique is a complete subgraph of any graph $G$. Clique complex is a simplicial complex formed from a network with the vertex set of $G$. If $c$ be the any clique in the network with $k$ vertices, then it will be (*k-1*)-simplex in the clique complex.

***Definition 15:*** Let $G = (V, E)$ be the underlying graph. A clique complex $X(G)$ is an abstract simplicial complex formed by the set of vertices $V$ of $G$, whose faces are all complete subgraphs of $G$.

It means a clique complex $X(G)$ is an abstract simplicial complex which follows the downward closure property. Vertex set V of graph G become the vertex set for clique complex $X(G)$ and every clique of size k become a simplex of dimension $(k-1)$. As $X(G)$ is an abstract simplicial complex it must follow the given conditions.

1) If $\sigma_i, \sigma_j \in X(G)$ then, either $\sigma_j \cap \sigma_j \in X(G)$ or $\sigma_j \cap \sigma_j = \emptyset$
2) $\forall \sigma_i \in X(G), (\sigma_k \subset \sigma_i \Rightarrow \sigma_k \in X(G))$. That is, every face of a clique $\sigma_i$ in $X(G)$ must also belong to $X(G)$.

The significance of the first condition is that, on arbitrarily choosing two cliques from $X(G)$, either they intersect at common face of the $X(G)$ or their intersection is $\emptyset$. The second condition signifies that for every arbitrary clique $c$ of $X(G)$, is closed under taking subsets. By the Definition 15 the construction of clique complexes of carbon mediation networks can be done. For example, the River-1 network from the aquatic environment it has 18 vertices and 33 edges. As in the clique complex all faces represents the clique of corresponding network, so the cliques of all sizes have been calculated. This network has 33 cliques of size 2 and 5 cliques of size 3.

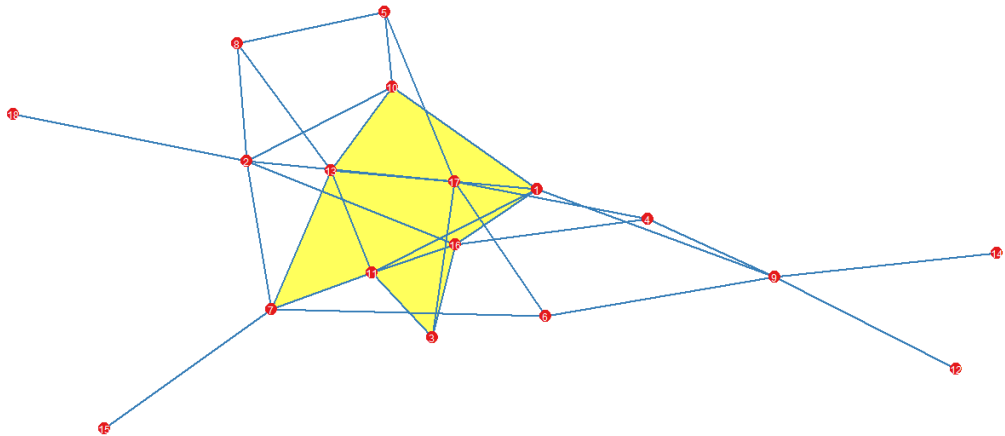

Clique complex for River-1 network

Clique complex $X(R)$ contains eighteen 0-simplcies, thirty-three 1-simplices which are the edges $\{1,10\}, \{1,9\}, \{1,11\}, \{1,13\}$ etc and five 2-simplices $\{7,11,13\}, \{3,11,16\}, \{1,11,16\}, \{1,11,13\}, \{1,10,13\}$. However, every $k$-clique is a $(k-1)$-simplex in the corresponding clique complex, which implies every edge will be a 1-simplex and every triangle will be a 2-simplex in $X(R)$.

We will look at non-pairwise interactions up to triangles in this work, with a concentration on interactions up to 3-cliques or 2-simplices. Let the clique complex of River-1 is $X(R)$, and the vertex set for this clique complex will be the same as vertex set of the River-1 network.

$X(R) = \{\{1\}, \{2\}, \dots, \{18\}, \{7,11\}, \{1,10\}, \dots, \{7,11,13\}, \{3,11,16\}, \dots\} = <1,2,3,\dots,32>$

$X(R)$ consists all 0-simplices, 1-simplices and 2-simplices from the River-1 network. It can be seen that the $X(R)$ satisfies both the conditions 1&2, it is closed under taking subsets. As we know that intersection of any two 0-simplices from $X(R)$ will be $\emptyset$, intersection of two 1-simplicies either 0-simplex $\emptyset$ and intersection of two 2-simplices either a 1-simplex, 0-simplex or $\emptyset$.

On arbitrarily choosing any two simplices of $X(R)$, their intersection i.e., $\{7,11,13\} \cap \{3,11,16\} = \{11\} \in X(R)$, which is again a face of $X(R)$. Hence, $X(R)$ is a Clique complex. In a similar way the clique complex of all the 10 networks can be constructed.

On the one hand, graph theory makes it clear when two vertices are adjacent to each other. Adjacency, on the other hand, is more difficult to define in simplicial complexes. Using the upper and lower adjacency criteria, it can be assessed whether the provided two $k$-simplices are adjacent or not. We will look at an example to better grasp the simplicial adjacency. First, the adjacency matrix at the 1-simplex level, i.e., between the simplicial complex's edges, will be computed.

This clique complex contains thirty-three 1-simplices, implying that the adjacency matrix at the 1-simplex level will be a square matrix of dimension thirty-three. Now, for the adjacency of these 1-simplices we need to check the upper adjacency as well as lower adjacency. For example, let's take two 1-simplices form river-1 network, $e_1 = \{7,15\}$ and $e_2 = \{7,13\}$ these 1-simplices are lower adjacent to each other because they have a common 0-simplex(vertex) or have a $(k-1)$-simplex in common (here $k=1$), hence it may be concluded that the $e_1$ and $e_2$ are lower adjacent to each other. However, for checking the upper adjacency of $e_1$ and $e_2$, by the Definition 10 in simplicial complex it can concluded that 1-simpices $\{7,15\}, \{7,13\}$ are upper adjacent to each other, because they are not a common face of a 2-simplex which is $\{7,13,15\}$. Because for two $k$-simplices to be upper adjacent they should not be the part of a $(k+1)$-simplex (here $k=1$). Hence 1-simplices $\{7,15\}, \{7,13\}$ are lower adjacent to each other but not upper adjacent to each other, which implies $\{7,15\}, \{7,13\}$ are adjacent to each other.

### 3.2 Simplicial centralities measures

We must first understand how the shortest path in a simplicial complex is defined before we can calculate the simplicial centralities. Because the centrality measures like closeness and betweenness depends on the shortest path between two $k$-simplices. As it is known previously that in graph theory shortest path between two vertices $u$ and $v$ is defined as the length of the shortest walk between $u$ and $v$. Walk between two vertices $u$ and $v$ is given by sequence of vertices and edges as $u, e_1, u_1, e_2, u_2, \dots, u_{n-1}, e_n, v$ where $n \in \mathbb{N}$. (2)

A few points must be understood before generalising the concept of walk for simplicial complexes.
- The sequence of vertices and edges in the definition of walk, in graph theory can be seen as sequence of 0-simplices and 1-simplices.

- Adjacency in graph theoretic sense is upper adjacency in simplicial complex, as defining a lower adjacency in graph theory would require the concept of a negative dimension.
- Two alternate consecutive 0-simplex are upper adjacent to each other in (2). (As we cannot define the lower adjacency for 0-simplex).

On the basis of above-mentioned points, the walk between two $k$-simplices $\sigma_i$ and $\sigma_j$ can be defined.

***Definition 16:*** A walk $W^k$ between two $k$-simplices ($k \geq 1$) $\sigma_i$ and $\sigma_j$ is the alternating sequence of $k$-simplices and $(k-1)$ simplices as $\sigma_i = \sigma_1, \alpha_1, \sigma_2, \ldots, \sigma_{n-1}, \alpha_{n-1}, \sigma_n = \sigma_j$ where each $\alpha_l$ is a $(k-1)$-simplex for $l = \{1,2,3, \ldots, n-1\}$, lower adjacent to $\sigma_l$ and $\sigma_{l+1}$, and two alternating consecutive $k$-simplices $\sigma_l$ and $\sigma_{l+1}$ cannot be upper adjacent to each other for $l = 1, \ldots, n-1$.

For $k = 0$, a walk on the 0-simplices (vertices) is a walk which has been defined in (2), and is the graph theoretic walk.

***Definition 17:*** Shortest path between two $k$-simplices ($k \geq 1$) $\sigma_i$ and $\sigma_j$ is a $W^k$ walk $\sigma_i, \alpha_1, \sigma_2, \ldots, \sigma_{n-1}, \alpha_{n-1}, \sigma_j$. Where $n$ is the minimum of all possible such sequences (that is, the sequence where each simplex appears only once). The value $n$ is the shortest path length between $\sigma_i$ and $\sigma_j$ is denoted by $d(\sigma_i, \sigma_j) = n$.

***Theorem 2:*** Shortest path length $d$ between two $k$-simplices $\sigma_i$ and $\sigma_j$ is a metric and $(\Delta, d)$ is a metric space.
***Proof:*** First we will prove that $d$ is metric on the elements of $\Delta$.
Let $\Delta$ be simplicial complex, a metric on $\Delta$ is a map $d: \Delta \times \Delta \rightarrow \mathbb{R}$
For $\forall \sigma_i, \sigma_j \in \Delta$

i. $d(\sigma_i, \sigma_j) \geq 0$ as the walk between them $\sigma_1, \alpha_1, \sigma_2, \ldots, \sigma_{n-1}, \alpha_{n-1}, \sigma_n$ always contains $n \geq 0$ $\alpha_i$s.
ii. $d(\sigma_i, \sigma_j) = d(\sigma_j, \sigma_i)$, as shortest path length between $\sigma_i$ and $\sigma_j$ is a walk $\sigma_i, \alpha_1, \sigma_2, \ldots, \sigma_{n-1}, \alpha_{n-1}, \sigma_j$ which can be written as $\sigma_j, \alpha_1, \sigma_2, \ldots, \sigma_{n-1}, \alpha_{n-1}, \sigma_i$ this walk represents the shortest path length of $\sigma_j$ and $\sigma_i$. Hence symmetry follows.
iii. For $\forall \sigma_i, \sigma_j, \sigma_k \in \Delta$, $d(\sigma_i, \sigma_j) + d(\sigma_j, \sigma_k) \geq d(\sigma_i, \sigma_k)$
Let $d(\sigma_i, \sigma_j) = n$, $d(\sigma_j, \sigma_k) = m$
$d(\sigma_i, \sigma_j) = n$, it represents the walk of length $n$ and $d(\sigma_j, \sigma_k) = m$, represents the walk of length $m$. The one minimum length walk between $\sigma_i, \sigma_k$ will be $\sigma_i, \alpha_1, \sigma_2, \ldots, \sigma_n, \alpha_n, \sigma_j, \beta_1, \ldots, \sigma_m, \beta_m, \sigma_k$ which has the length $(n + m)$. So, no other walk of length greater than $(n + m)$ can be the shortest path length walk for $\sigma_i$ and $\sigma_k$, the only possible minimum length walk between $\sigma_i$ and $\sigma_k$, can be of length less than $(n + m)$. Hence triangle inequality follows.

Hence $d$ follows the all three conditions of a metric $\Rightarrow (\Delta, d)$ is a metric space. ∎

For example, according to Definition 16, a walk between two 1-simplices $\sigma_i$ and $\sigma_j$ from Lake-1 network is an alternating sequence of 1-simplices and 0-simplices such that it where $\sigma_i = \{15,10\}, \sigma_j = \{8,16\}$. The walk $W^1$ may be written as
$W^1 = \{15,10\}; \{10\}, \{10,6\}; \{6\}; \{6,8\}; \{8\}; \{8,16\}$.
Note that there exists more than one walk between any two simplices in this example, and hence we can write the $W^2$ walk between two 2-simplices $\rho = \{11,14,12\}$ and $\varphi = \{12,14,1\}$ as $W^2 = \{11,14,12\}; \{14,12\}, \{12,14,8\}; \{12,14\}; \{14,12,1\}$.

In the walk $W^2$, $\{14,12\}$ is a 1-simplex which is a common face for both the simplices $\{11,14,12\}$ and $\{12,14,8\}$, and these 2-simplices are not upper adjacent to each other. Similar argument is valid for next terms in the sequence $W^2$. Hence $W^2$ is the minimum length walk between $\rho$ and $\varphi$.

Connectivity in simplicial complex can be defined at various levels designated as $k$, where $k = 0,1,2,...$. For $k = 0$ it represents the graph-theoretic connectivity. A graph is said to be connected if for every pair of vertices there exist a path. Similarly, the concept of connectedness in simplicial complexes can be generalised.

***Definition 18:*** Simplicial complex $\Delta$ is defined as connected iff for every two k-simplices $\sigma_i$ and $\sigma_j$ there exist a walk $\sigma_i, \alpha_1, \sigma_2, ..., \sigma_{n-1}, \alpha_{n-1}, \sigma_j$ of finite length, such that $d(\sigma_i, \sigma_j) =$ finite.

***Theorem 3:*** Adjacency matrix $A^k$ is irreducible if and only if $\Delta$ is a connected simplicial complex at $k^{th}$ level.
***Proof:*** A matrix of order $n$ is said to be reducible if and only if there are two disjoint sets of indexes $I_1$ and $I_2$ such that there exist square matrix blocks of size equal to the cardinality of $I_1$ and $I_2$ and $|I_1|+|I_2|=n$ such that for every $(i,j) \in I_1 \times I_2$, $a_{ij} = 0$. Otherwise, it said to be irreducible.
We are given that $A^k$ is irreducible $\Rightarrow \neg \exists$ any block partition with respect to index sets $I_1$ and $I_2$ such that for every $(i,j) \in I_1 \times I_2$, $a_{ij} = 0$. This implies that for any two given k-simplices $\sigma_i$ and $\sigma_j$ in $\Delta$ either $a_{ij} \neq 0$ or $a_{ij} = 0$. If $a_{ij} \neq 0$ then $\sigma_i$ and $\sigma_j$ are adjacent to each other hence there is a walk of length one between $\sigma_i$ and $\sigma_j$. If $a_{ij} = 0$ then $\exists m \in \mathbb{N}$ such that the $a_{ij}^{th}$ entry of $A^k \times A^k \times ... \times A^k$ (m times) is non-zero (as the matrix $A^k$ is symmetric non-negative matrix) [18]. Because $A^k$ is an irreducible matrix it gives us assurance for the existence of such $m \in \mathbb{N}$ otherwise $\sigma_i$ or $\sigma_j$ become the isolated k-simplex which gives a zero row/column in the irreducible matrix $A^k$, which is a contradiction. It means there exist a walk of length $m$ between $\sigma_i$ and $\sigma_j$. Hence $\Delta$ is connected at $k^{th}$ level.
Conversely, let us assume $A^k$ is a reducible matrix of order n. This implies there exist two square matrix blocks with respect to the index sets $I_1$ and $I_2$ such that $|I_1| + |I_2| = n$. Let $I_1 = \{1,2,...,m\}$ and $I_2 = \{m+1, m+2, ..., m+l\}$ and $m + l = n \in \mathbb{N}$.

Now for $(1, m + 1) \in I_1 \times I_2$ the entry $a_{1(m+1)} = 0$ which implies that $\exists$ two k-simplices $\sigma_1$ and $\sigma_{1+m}$ which are not adjacent to each other.

However, if a matrix $A$ is reducible then $A \times A \times ... \times A$ ($p$ times) is reducible for any $p \in \mathbb{N}$. It means $A^k \times A^k \times ... \times A^k$ ($p$ times) is reducible for every $p \in \mathbb{N}$. Hence there does not exist any walk of finite length between $\sigma_1$ and $\sigma_{1+m}$. But we are given that $\Delta$ is a connected simplicial complex at $k^{th}$ level, which is a contradiction to our assumption.

Hence adjacency matrix $A^k$ of simplicial complex $\Delta$ is irreducible. ∎

***Corollary 1:*** Clique complex $X(G)$ of a connected graph $G$ need not to be connected.

We shall illustrate this corollary through an example.

let $G$ be a connected graph and $X(G)$ be its clique complex.

We can clearly see that $G$ is connected as we can find a path of finite length between any two vertices of $G$.

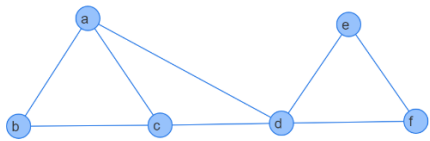

$$X(G) = \left\{ \begin{array}{c} \{a\}, \{b\}, \{c\}, \{d\}, \{e\}, \{ab\}, \{ac\} = e_1, \{ad\}, \{bc\}, \{cd\}, \{df\} = e_2, \{fe\}, \\ \{de\}, \{abc\}, \{acd\}, \{def\} \end{array} \right\}$$

$X(G)$ is the clique complex of $G$.

Claim: $X(G)$ is not connected at k=1 level, as there does not exist any walk of finite length between 1-simplex $e_1$ and $e_2$. Moreover, there does not exist any walk of finite length from $e_1$ to every other 1-simplex of $X(G)$. It implies 1-simplex $e_1$ is an isolated simplex and row corresponding to $e_1$ in adjacency matrix $A^1$ is a zero row. Which implies $A^1$ cannot be an irreducible matrix.

Hence by Theorem 3, $X(G)$ is a disconnected clique complex at k=1 level. ∎

In the following paragraphs, we discuss some of the important centrality indices that helps us to obtain a deeper understanding of the structure of the ecological networks that we construct, with the results of construction being detailed in Section 4.

**Degree centrality**

The most common centrality measure for a network is the degree centrality. In graphs degree of a vertex $v$ defined as the number of adjacent vertices to $v$ but for a simplicial complex we have certain levels of degree measures, that can be designated as $\delta_k(i)$, where $k = \{0,1,2,...\}$ is the level of the simplex and $i$ is the corresponding simplex, for example $k = 0$ represents the degree at vertex level or 0-simplex level, k=1 represents the degree at 1-simplex level, k=2 represents the degree at 2-simplex level. The degree of a $k$-simplex is the number of other k-simplices to which it is adjacent. So, for calculating the degree centrality at different levels, adjacency matrices $A^k$ at different levels need to be calculated. Where

k=0, it represents the adjacency matrix of a graph.

***Theorem 4:*** Degree of a $k$-simplex $\sigma_i$ in a simplicial complex $\Delta$ is given by the corresponding row sum in the adjacency matrix $A^k$.

***Proof:*** $A^k$ represents the adjacency matrix of $k$-simplices in $\Delta$, it is a non-negative matrix with entries either 0 or 1. If two $k$-simplices are adjacent to each other it gives the value 1 otherwise 0. For a $k$-simplex $\sigma_i$, the row corresponding to $\sigma_i$ tells us about the adjacency and non-adjacency of $\sigma_i$, entries with 0 value tells us about the non-adjacency of $\sigma_i$ and entries with value 1 tells about the adjacency of $\sigma_i$, if we take the row sum corresponding to $\sigma_i$ which will be the sum of $0's$ and $1's$ can tell us about the degree of $\sigma_i$ because it will let us know the number $k$-simplices to which $\sigma_i$ is adjacent. ∎

**Closeness centrality**
Every centrality measure gives different information about the vertices of any graph. For example, closeness centrality informs about those vertices which are responsible for spreading information in the given graph. This concept of closeness centrality can be generalised for the simplicial complex obtained in our modelling.

The graph-theoretic closeness centrality depends upon shortest path distance between the vertices, and measures the closeness of a given vertex in terms of the shortest path distance from other vertices of the network. A generalization of this concept for simplicial structures analogously must measure the closeness of a given $k$-simplex in terms of shortest path distance with other $k$-simplices. $k$-simplices with a high closeness score have the shortest distances to all other $k$-simplices. The definition of simplicial closeness depends on the concept of simplicial farness, and cannot be defined for disconnected simplicial complexes. The simplicial farness depends upon shortest path distance of the given $k$-simplex to other $k$-simplices in the simplicial complex. So, it can be defined as, for any given $k$–simplex $\sigma_i$ the simplicial farness is the sum of its shortest path distances to all other $k$-simplices in the simplicial complex. Simplicial farness for $\sigma_i$ is given as

$$F(\sigma_i) = \sum_{\sigma_i \neq \sigma_j} d(\sigma_i, \sigma_j)$$

The simplicial closeness is the reciprocal of simplicial farness. Simplicial closeness $C$ for $k$-simplex $\sigma_i$ given by

$$C(\sigma_i) = \frac{1.}{\sum_{\sigma_i \neq \sigma_j} d(\sigma_i, \sigma_j)}.$$

Where $d(\sigma_i, \sigma_j)$ is the shortest path distance between $\sigma_i$ and $\sigma_j$.

**Betweenness centrality**
Betweenness centrality also depends upon the shortest path distance between simplices. It addresses the k-simplices that influence the simplicial complex, as well as the k-simplices that control the flow of information in the simplicial complex. The expression for simplicial betweenness centrality can be directly generalised from the graph-theoretic expression.

The betweenness centrality of a $k$-simplex $\sigma_j$ is given by the expression:

$$B(\sigma_j) = \sum \frac{N_{(\sigma_i,\sigma_k)}(\sigma_j)}{N_{(\sigma_i,\sigma_k)}}$$

Where $N_{(\sigma_i,\sigma_k)}$ is the total number of shortest paths between $\sigma_i$ and $\sigma_k$ and $N_{(\sigma_i,\sigma_k)}(\sigma_j)$ is the total number of shortest paths between $\sigma_i$ and $\sigma_k$ through $\sigma_j$.

**Eigenvector centrality**
Eigenvector centrality measures a vertex's importance while giving consideration to the importance of its neighbors. A vertex with very few connections can be having very high eigenvector centrality if those few are very well connected to other vertices in the network. In the graph theory eigenvector score of any $i^{th}$ vertex $v_i$ of any graph $G$ is given by the $i^{th}$ component of the principal eigenvector of the adjacency matrix $A$.
The relative centrality score of a vertex $v$ is given by

$$x_v = \frac{1}{\lambda} \sum_{u \in N(v)} x_u$$

Where, $\lambda$ is a scalar and $N(v)$ is the set of neighbors of vertex $v$.
The simplicial eigenvector centrality of any $k$-simplex $\sigma_i$ by the $i^{th}$ component of the principal eigenvector of the adjacency matrix $A^k$.

***Definition 19:*** A subcomplex $L$ of a simplicial complex $\Delta$ is the simplicial complex such that $L \subseteq \Delta$.

**Subgraph centrality**
Subgraph centrality measures a vertex's importance on the basis of the ability to participate in the subgraphs of the network. A generalization of this concept for simplicial structures analogously must measure the importance of $k$-simplex on the basis of the ability to participate in the in the subcomplex of the simplicial complex at k[th] level [64, 65].
Subgraph centrality of vertex $v$ in any network G is given by

$$SC(v) = \sum_{m=0}^{\infty} \frac{\mu_m(v)}{m!}$$

where,
   $\mu_m(v)$ = Number of closed walks of length $m$ starting and ending at same vertex $v$.
   $\mu_m(v) = [A^m]_{ii}$, $[A^m]_{ii}$ is the i[th] diagonal entry of m[th] power of adjacency matrix $A$.

Similarly, for a simplicial complex $\Delta$ with adjacency matrix $A^k$ at k[th] level, the subgraph centrality for a $k$-simplex $\sigma$ is given by

$$SC(\sigma) = \sum_{m=0}^{\infty} \frac{\mu_m(\sigma)}{m!}$$

where,

$\mu_m(\sigma)$ = Number of closed walks of length $m$ starting and ending at same $k$-simplex $\sigma$.

$\mu_m(\sigma) = [A^k]^m{}_{ii}$, $[A^k]^m{}_{ii}$ is the $i^{th}$ diagonal entry of $m^{th}$ power of adjacency matrix $A^k$.

## 4 Results

1) A summary of the interaction networks originating from the aquatic environment and terrestrial environment studied in this work presented in Table1.

2) Except for the networks of Lake-3 and Forest-2, all the clique complexes corresponding to each network analysed in this paper are connected at the 1-simplex level, according to Theorem 3. However, their corresponding clique complexes show different connectivity at different $k$-simplex levels (Table19).

3) Table 3 and 4 show that for the centralities at 0-simplex level of the Caribbean Food Web network vertex $v_1$ has the highest score for all centrality indices, which concluded that $v_1$ is the specie/taxa which involved in the most number of carbon mediations with other ecological units. Thus, it may be inferred that $v_1$ is the most important vertex 0-simplex level of the network. Further, after constructing the corresponding clique complex at 1-simplex and 2-simplex levels, it can be observed that the top ranked $k$-simplices are those that are involved in grouping with $v_1$. Hence, it can be concluded that $v_1$ participated in lower order interaction as well as in the higher-order interaction simultaneously. Consequently $v_1$ remains a very highly central vertex in both, the 1-skeleton as also the corresponding clique complex.

4) Table 3 and 4 show that the low centrality scores of some 0-simplices corresponding to top ranked 1-simplices and 2-simplices shows less involvement in the mediation of carbon at the 0-simplex level. By calculating the centrality scores at different levels, it may be concluded that the bottom ranked 0-simplices can also be members of top ranked simplices at different levels. As a result, certain ecological units which are not significantly involved (that is, have lower rank) in the pairwise interaction can participate in higher order interactions.
Top ranked 1 and 2-simplices of Caribbean Food Web network are found to be:

$e_{16}$ is the 1-simpex of $v_1$ and $v_{26}$
$e_{21}$ is the 1-simplex of $v_1$ and $v_{32}$
$e_{13}$ is the 1-simplex of $v_{23}$ and $v_1$
$e_{82}$ is the 1-simplex of $v_7$ and $v_{14}$

$t_1$ is 2-simplex of vertex $v_1$, $v_{29}$ and $v_4$
$t_{13}$ is 2-simplex of $v_1$, $v_4$ and $v_{42}$
$t_{49}$ is 2-simplex vertex $v_{15}$, $v_{31}$ and $v_{36}$

The ecological units corresponding to vertices $v_i$ are as given below:

$v_1$: Anolis gingivinus
$v_7$ : Yellow warbier
$v_{14}$: Coleoptera adult
$v_{26}$: Annelid
$v_{23}$: Other hymenoptera
$v_{32}$: Diptera adult
$v_4$ : Anolis pogus

$v_{29}$ : Millipede
$v_{36}$ : Leaves
$v_{42}$ : Thalandros cubensis
$v_{15}$ : Orthoptera
$v_{31}$ : Homoptera

5) Fig.7(b) represents the clique complex of Caribbean Food-Web network, in which two 1-simplices are adjacent if they share a common 0-simplex and not part of a 2-simplex in network. This implies that two 1-simplices are interacting at 1-simplex level if they have a common carbon interactant and are not involved in 2-simplex of the carbon mediation network. This result shows the significance of the adjacency of 1-simplices in the network.

6) It may be observed from Tables 9 & 10 that every top ranked 1-simplices $\{v_2, v_{11}\}$, $\{v_2, v_9\}$, $\{v_2, v_{12}\}$, $\{v_{12}, v_{19}\}$ and $\{v_{13}, v_{19}\}$ from the Lake-1 network contains top ranked 0-simplices $\{v_2\}$ and $\{v_{19}\}$. Also, the 0-simplices other than top ranked $\{v_{11}\}$, $\{v_{12}\}, \{v_{13}\}, \{v_9\}$ are not participating in the important higher order interactions. The top-ranked 2-simplices which are involved in higher order interactions contains middle-ranked 0-simplices at 0-simplex level.

7) It may be observed from Tables 9 and 10, $\{v_2\}$ is the top ranked 0-simplex and is also involved with the top ranked 1-simplices but it is not participating in the top ranked higher order interactions. Thus, it may be concluded that higher order interactions that include $\{v_2\}$ are not of importance from carbon mediation point of view. Further, with the help of this centrality score those ecological units which have the maximum involvement in higher order interactions can be identified. These are:

$v_2$: Cryptomonas sp. 3 - Cryptomonas sp. 4 - Cosmarium sp. - Dactylococcopsis fascicularis -Dictyosphaerium pulchellum - Dinobryon sertularia - Sphaerocystis schroeteri - Glenodinium
pulvisculus - Oocystis sp. 1 - Oocystis sp. 2 - Peridinium pulsillum - Schroederia setigera
$v_9$: Conochiloides colonial
$v_{12}$: Daphnia pulex
$v_{13}$: Daphnia rosea
$v_{11}$: Diaptomus oregonensis

$v_{19}$: *Chaoborus punctipennis*

Similar results are obtained for all 10 carbon mediation networks.

8) By plotting the graph centralities and simplicial centralities for all the networks studied in this work, it may be concluded that the ranking of 0-simplices, 1-simplices and 2-simplices according to a given centrality measure can differ significantly for every simplicial complex, and observed to be independent of the type (terrestrial or aquatic) of the corresponding networks. Two plots as representative sample are given as Plot 1 and Plot 2, for River-1 and Caribbean Food Web network respectively.

9) The simplicial centrality plot for every network shows that top ranked 0-simplices appears in top ranked 1-simplices and 2-simplices. It therefore, may be concluded that 0-simplices with high centrality scores are likely to participate in the higher order interactions, with the exception of some ecological units.

10) From Table 20, it may be observed that Lake-3 network is a disconnected network, and hence, "Dissolved organic carbon" is an isolated 0-simplex. Corresponding clique complex is also disconnected at 1-simplex and 2-simplex level. 1-simplex $e_1$ between $v_1$ and $v_5$ (*Small ciadocerans and Copepods*) is an isolated 1-simplex. By using Theorem 3, it may be concluded that adjacency matrices $A^0$, $A^1$ and $A^2$ are reducible matrices. It can be identified the 0-simplices of $e_1$ as the ecological units which involved in pair wise interaction only because they cannot form any triangle(2-simplex) which represents the higher order interaction in the clique complex.

11) Table 20 further shows that the Forest-2 network is a disconnected network, but the corresponding clique complex is connected at 1-simplex level and disconnected at 2-simplex level. Hence, it may be inferred that for terrestrial ecosystems, existence of simplicial connectivity at one level does not imply the existence of connectivity at other levels.

**Table 1** Summary of order, size of carbon mediation networks.

| S.No. | Network | Environment | No. of vertices | No. of edges | Reference |
|---|---|---|---|---|---|
| 1 | Estuary-2 | Aquatic | 25 | 77 | 66 |
| 2 | Lake-1 | Aquatic | 20 | 51 | 67 |
| 3 | Lake-3 | Aquatic | 21 | 66 | 68 |

| 4 | River-1 | Aquatic | 18 | 33 | 69 |
| 5 | Marine-1 | Aquatic | 28 | 194 | 70 |
| 6 | Rocky shore-1 | Terrestrial | 21 | 28 | 71 |
| 7 | Rocky shore-4 | Terrestrial | 20 | 55 | 72 |
| 8 | Forest-2 | Terrestrial | 31 | 58 | 73 |
| 9 | Forest-4 | Terrestrial | 29 | 65 | 74 |
| 10 | CFW | Terrestrial | 43 | 207 | 75 |

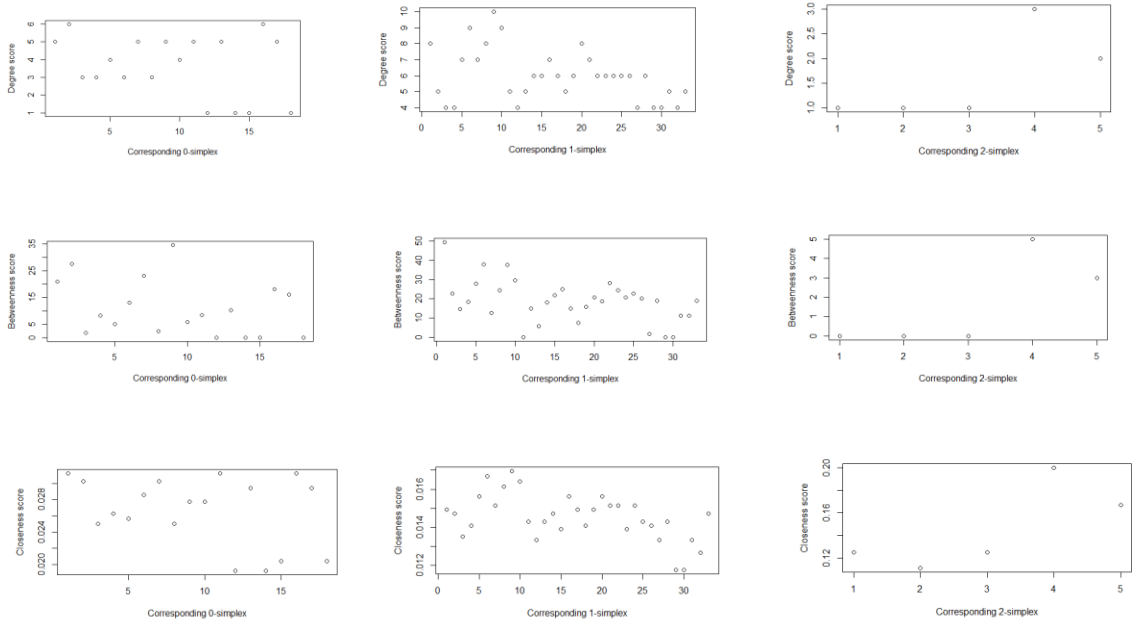

**Plot (1):** Comparison of rankings of 0-simplices, 1-simplices and 2-simplices for River-1 network

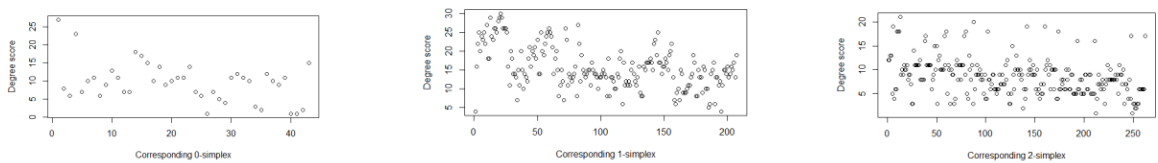

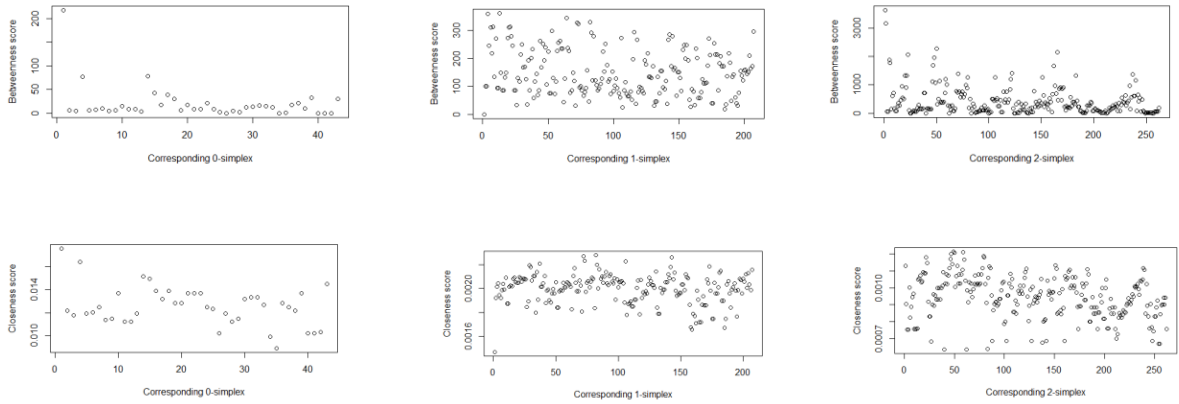

**Plot (2):** Various centralities plot of 0-simplices, 1-simplices and 2-simplices of Caribbean Food web network for aquatic environment

**Table 3: Caribbean food-web network from terrestrial environment**

| Network | Number of vertices | Number of edges | Number of 2 cliques | Number of 3 cliques | Clique number |
|---|---|---|---|---|---|
| Caribbean Food web | 43 | 207 | 207 | 262 | 6 |

**Table 4: Comparison of centralities at all three simplex levels of Caribbean Food Web network**

| Network | Highest subgraph centrality and simplex | Highest degree centrality and simplex | Highest eigenvector centrality and simplex | Highest betweenness centrality and simplex | Highest closeness centrality and simplex |
|---|---|---|---|---|---|
| 0-simplex level | 23905.80, $v_1$ | 27, $v_1$ | 1, $v_1$ | 217.204, $v_1$ | 0.017, $v_1$ |
| 1-simplex level | 106489716, $e_{16}$ | 30, $e_{21}$ | 1.0, $e_{16}$ | 360.87, $e_{13}$ | 0.00227, $e_{82}$ |
| 2-simplex level | 1118011, $t_{13}$ | 21, $t_{13}$ | 1.0, $t_1$ | 3629.91, $t_1$ | 0.00121, $t_{49}$ |

$e_{16}$ is the 1-simplex of $v_1$ and $v_{26}$
$e_{21}$ is the 1-simplex of $v_1$ and $v_{32}$
$e_{13}$ is the 1-simplex of $v_{23}$ and $v_1$
$e_{82}$ is the 1-simplex of $v_7$ and $v_{14}$
$t_1$ is 2-simplex of vertex $v_1$, $v_{29}$ and $v_4$

$t_{13}$ is 2-simplex of $v_1, v_4$ and $v_{42}$
$t_{49}$ is 2-simplex vertex $v_{15}, v_{31}$ and $v_{36}$

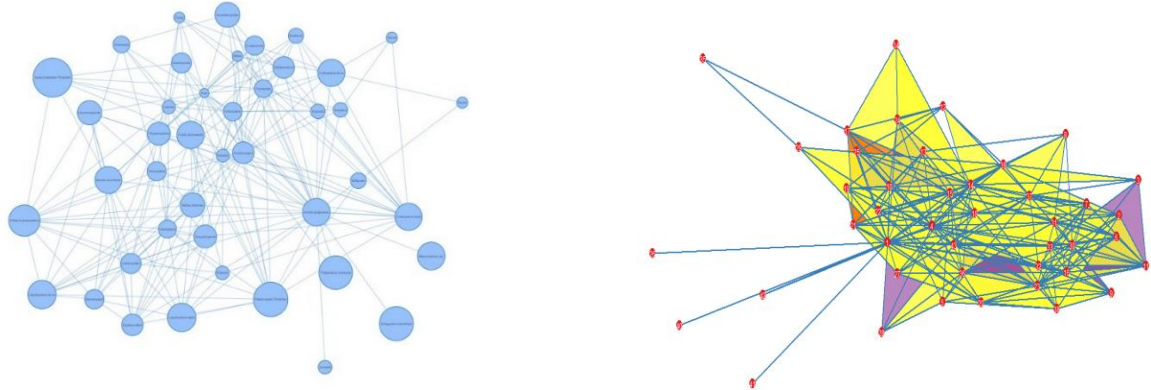

(a) Pairwise network of Caribbean Food Web     (b) Clique complex of Caribbean Food Web

Fig.7: Pairwise network and corresponding clique complex

**Table 5: River-1 network from aquatic environment**

| Network | Number of vertices | Number of edges | Number of 2 cliques | Number of 3 cliques | Clique number |
|---|---|---|---|---|---|
| River 1 | 18 | 33 | 33 | 5 | 3 |

**Table 6: Comparison of centralities at all three simplex levels of River-1 network**

| Network | Highest subgraph centrality and simplex | Highest degree centrality and simplex | Highest eigenvector centrality and simplex | Highest betweenness centrality and simplex | Highest closeness centrality and simplex |
|---|---|---|---|---|---|
| 0-simplex level | 12.77, $v_{16}$ | 6, $v_{16}$ | 1, $v_{16}$ | 34.57, $v_9$ | 0.03125, $v_{16}$ |
| 1-simplex level | 98.19, $e_9$ | 10, $e_9$ | 1.0, $e_9$ | 49.43, $e_1$ | 0.0169, $e_9$ |
| 2-simplex level | 2.9668, $t_4$ | 3, $t_4$ | 1, $t_4$ | 5, $t_4$ | 0.20, $t_4$ |

$e_9$ is 1-simplex of $v_2$ and $v_{16}$
$e_1$ is 1-simplex of vertex $v_1$ and $v_9$
$t_4$ is 2-simplex of $v_1$, $v_{11}$, $v_{13}$

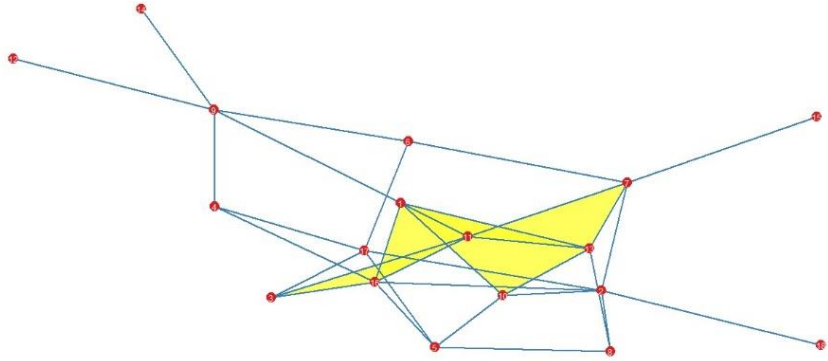

Fig.8: Clique complex of River-1 network

**Table 7: Estuary-2 network from aquatic environment**

| Network | Number of vertices | Number of edges | Number of 2 cliques | Number of 3 cliques | Clique number |
|---|---|---|---|---|---|
| Estuary- 2 | 25 | 77 | 77 | 53 | 4 |

**Table 8: Comparison of centralities at all three simplex levels of Estuary-2 network**

| Network | Highest subgraph centrality and simplex | Highest Degree centrality and simplex | Highest eigenvector centrality and simplex | Highest betweenness centrality and simplex | Highest closeness centrality and simplex |
|---|---|---|---|---|---|
| 0-simplex level | 537.68, $v_{25}$ | 21, $v_{25}$ | 1, $v_{25}$ | 131.211, $v_{25}$ | 0.037, $v_{25}$ |
| 1-simplex level | 2322461, $e_{23}$ | 21, $e_{77}$ | 1, $e_{23}$ | 225.921, $e_{70}$ | 0.0064, $e_{73}$ |
| 2-simplex level | 99.16, $t_{42}$, $t_{32}$ | 10, $t_{42}, t_{32}$ | 1, $t_{42}$ | 215.166, $t_{20}$, $t_{16}$ | 0.0084, $t_8$ |

$e_{23}$ is 1-simplex of $v_8$ and $v_{25}$
$e_{70}$ is 1-simplex of $v_{19}$ and $v_{24}$
$e_{73}$ is 1-simplex of $v_{20}$ and $v_{24}$
$e_{77}$ is 1-simplex of $v_{21}$ and $v_{25}$
$t_8$ is 2-simplex of $v_{14}$, $v_{18}$ and $v_{25}$
$t_{16}$ is 2-simplex of $v_{12}$, $v_{18}$ and $v_{25}$
$t_{20}$ is 2-simplex of $v_{11}$, $v_{18}$ and $v_{25}$

$t_{32}$ is 2-simplex of $v_{16}$, $v_{18}$ and $v_{25}$
$t_{42}$ is 2-simplex of $v_{15}$, $v_{18}$ and $v_{25}$

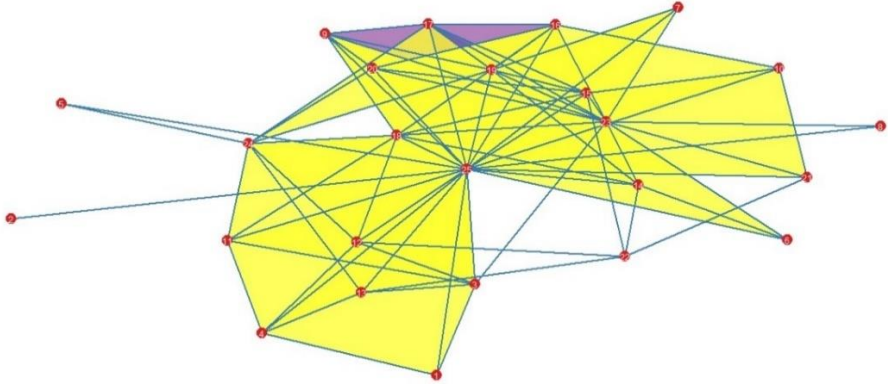

Fig.9: Clique complex of Estuary-2 network

**Table 9: Lake-1 network from aquatic environment**

| Network | Number of vertices | Number of edges | Number of 2 cliques | Number of 3 cliques | Clique number |
|---|---|---|---|---|---|
| Lake-1 | 20 | 51 | 51 | 24 | 4 |

**Table 10: Comparison of centralities at all three simplex levels of Lake-1 network**

| Network | Highest subgraph centrality and simplex | Highest degree centrality and simplex | Highest eigenvector centrality and simplex | Highest betweenness centrality and simplex | Highest closeness centrality and simplex |
|---|---|---|---|---|---|
| 0-simplex level | 102.99, $v_{19}$ | 11, $v_{19}$ | 1, $v_{19}$ | 40.29, $v_2$ | 0.03, $v_{19}$ |
| 1-simplex level | 1238.76, $e_6$ | 14, $e_{42}$ | 1, $e_6$ | 82.56, $e_{39}$ | 0.0107526, $e_6$ and $e_5$ |
| 2-simplex level | 33.88, $t_{12}$ | 7, $t_{12}$ | 1, $t_{12}$ | 56, $t_{11}$, $t_7$ | 0.0066, $t_{11}$ |

$e_6$ is 1-simplex of $v_2$ and $v_{11}$
$e_5$ is 1-simplex of $v_2$ and $v_9$
$e_7$ is 1-simplex of $v_2$ and $v_{12}$
$e_{39}$ is 1-simplex of $v_{12}$ and $v_{19}$
$e_{42}$ is 1-simplex of $v_{13}$ and $v_{19}$
$t_7$ is 2-simplex of vertex $v_{17}$, $v_{18}$ and $v_{19}$
$t_{11}$ is 2-simplex of $v_{12}$, $v_{17}$ and $v_{19}$
$t_{12}$ is 2-simplex vertex $v_{10}$, $v_{17}$ and $v_{19}$

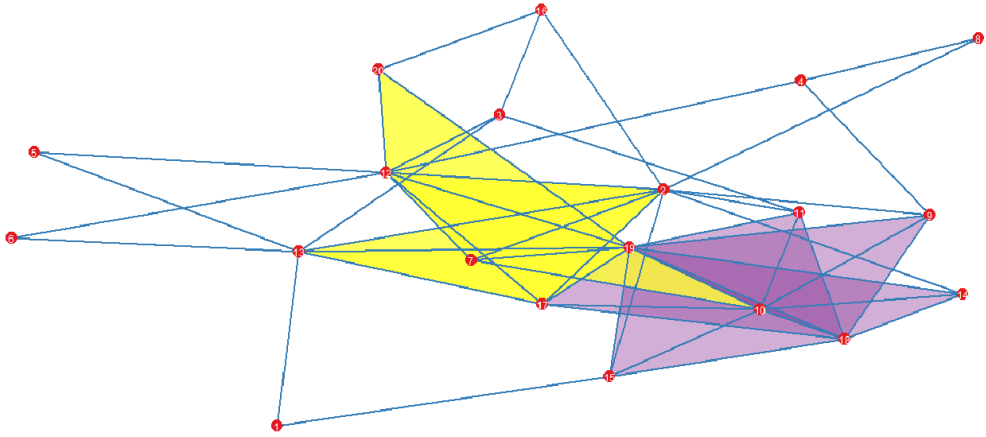

Fig.10: Clique complex of the Lake-1 network

**Table 11: Rocky shore-1 network from terrestrial environment**

| Network | Number of vertices | Number of edges | Number of 2 cliques | Number of 3 cliques | Clique number |
|---|---|---|---|---|---|
| Rocky shore-1 | 21 | 28 | 28 | 0 | 2 |

**Table12: Comparison of centralities at all three simplex levels of Rocky shore-1 network**

| Network | Highest subgraph centrality and simplex | Highest degree centrality and simplex | Highest eigenvector centrality and simplex | Highest betweenness centrality and simplex | Highest closeness centrality and simplex |
|---|---|---|---|---|---|
| 0-simplex level | 11.82, $v_4$ | 8, $v_4$ | 1, $v_4$ | 84.09, $v_1$ | 0.0217, $v_{18}$ |
| 1-simplex level | 208.901, $e_{17}$ | 10, $e_{17}$ | 1, $e_{17}$ | 89.833, $e_1$ | 0.01886, $e_1$ |

$e_1$= 1-simplex of $v_1$ and $v_2$
$e_{17}$= 1-simplex of $v_4$ and $v_{18}$

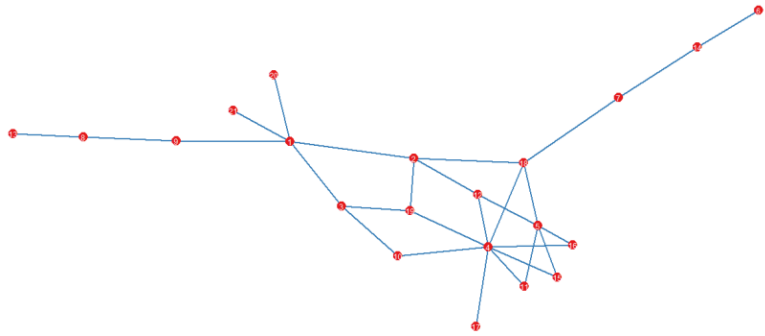

Fig.11: Clique complex of the Rocky Shore-1 network

**Table 13: Forest-4 network from terrestrial environment**

| Network | Number of vertices | Number of edges | Number of 2 cliques | Number of 3 cliques | Clique number |
|---|---|---|---|---|---|
| Forest-4 | 29 | 65 | 65 | 7 | 3 |

**Table14: Comparison of centralities at all three simplex levels of Forest-4 network**

| Network | Highest subgraph centrality and simplex | Highest degree centrality and simplex | Highest eigenvector centrality and simplex | Highest betweenness centrality and simplex | Highest closeness centrality and simplex |
|---|---|---|---|---|---|
| 0-simplex level | 50.110, $v_8$ | 10, $v_8$ | 1, $v_8$ | 79.64, $v_3$ | 0.01724, $v_8$ |
| 1-simplex level | 1241.55, $e_{22}$ | 13, $e_6$, $e_{22}$, $e_{27}$ | 1, $e_{22}$ | 122.44, $e_6$ | 0.0071, $e_6$ |
| 2-simplex level | 6.391, $t_6$ | 4, $t_6$ | 1, $t_6$ | 6, $t_6$ | 0.07, $t_6$ |

$e_{22}$ = 1-simplex of $v_8$ and $v_{22}$
$e_{27}$ = 1-simplex of $v_{21}$ and $v_9$
$e_6$ = 1-simplex of $v_3$ and $v_9$
$t_6$ = 2-simplex of $v_8$, $v_{18}$ and $v_{21}$

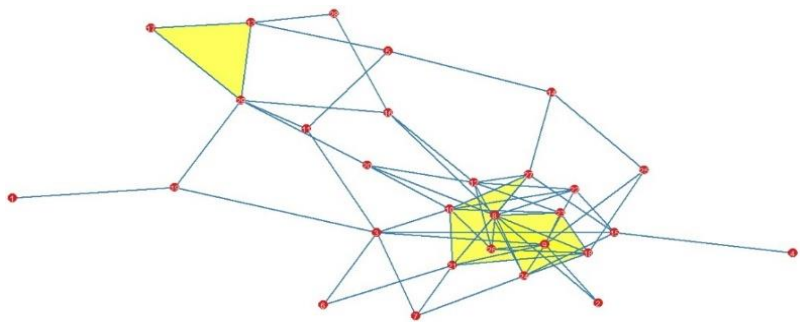

Fig.11: Clique complex of Forest-4 network.

**Table 17: Centralities at 0-simplex level of networks**

|  | Highest subgraph centrality | Highest degree centrality | Highest eigenvector centrality | Highest betweenness centrality | Highest closeness centrality |
|---|---|---|---|---|---|
| Forest-2 | 20.032 | 9 | 1 | 108.249 | 0.011 |
| Forest-4 | 50.110 | 10 | 1 | 79.64 | 0.01724 |
| Lake-3 | 281.311 | 11 | 1 | 20.63 | 0.20833 |
| Lake-1 | 102.99 | 11 | 1 | 40.29 | 0.03 |
| Caribbean Food Web network | 23905.80 | 27 | 1 | 217.205 | 0.017 |
| Estuary-2 | 537.68 | 21 | 1 | 131.211 | 0.037 |
| Rocky shore-1 | 11.82 | 11 | 1 | 84.09 | 0.0217 |
| Marine-1 | 510282.73 | 26 | 1 | 28.65 | 0.0333 |
| Rocky shore-4 | 137.44 | 13 | 1 | 43.09 | 0.04 |
| River-1 | 12.77 | 6 | 1 | 34.57 | 0.03125 |

**Table 18: Centralities at 1-simplex level of networks**

|  | Highest subgraph centrality | Highest degree centrality | Highest eigenvector centrality | Highest betweenness centrality | Highest closeness centrality |
|---|---|---|---|---|---|
| Forest-2 | 1413.25 | 15 | 1 | 151.90 | 0.0089 |
| Forest-4 | 1241.55 | 13 | 1 | 122.44 | 0.0071 |
| Lake-1 | 1238.76 | 13 | 1 | 82.56 | 0.0107526 |
| Lake-3 | 343.71 | 12 | 1 | 156.31 | 0.005 |
| Caribbean Food Web network | 106489716 | 30 | 1 | 360.87 | 0.0022 |
| Estuary-2 | 2322461 | 21 | 1 | 225.921 | 0.0064 |
| Rocky shore-1 | 208.901 | 10 | 1 | 89.833 | 0.01886 |
| Rocky shore-4 | 4997.09 | 15 | 1 | 68.58 | 0.009 |
| Marine-1 | 74048.62 | 23 | 1 | 992.82 | 0.00221 |
| River-1 | 98.19 | 10 | 1.0 | 49.43 | 0.0169 |

**Table 19: Centralities at 2-simplex level of networks**

|  | Highest subgraph centrality | Highest degree centrality | Highest eigenvector centrality | Highest betweenness centrality | Highest closeness centrality |
|---|---|---|---|---|---|
| Forest-2* | 0 | 0 | 0 | 0 | 0 |
| Forest-4 | 6.391, $t_6$ | 4, $t_6$ | 1, $t_6$ | 6, $t_6$ | 0.07, $t_6$ |
| Lake-1 | 33.88, $t_{12}$ | 7, $t_{12}$ | 1, $t_{12}$ | 56, $t_{11}, t_7$ | 0.0066, $t_{11}$ |
| Lake-3 | 32.89, $t_{48}$ | 9, $t_{48}$ | 1, $t_{48}$ | 435.544, $t_{33}$ | 0.0011, $t_{33}$ |
| Caribbean Food Web network | 1118011, $t_{13}$ | 21 | 1 | 3629.91 | 0.00121, $t_{49}$ |
| Estuary-2 | 99.16, $t_{42}, t_{32}$ | 10, $t_{42}, t_{32}$ | 1, $t_{42}$ | 215.166, $t_{20}, t_{16}$ | 0.0084, $t_8$ |
| Rocky shore-1** | - | - | - | - | - |
| Rocky shore-4 | 73.36, $t_{16}$ | 9, $t_{16}$ | 1, $t_{16}$ | 66, $t_2$ | 0.021, $t_{16}$ |
| Marine-1 | 38140, $t_{221}$ | 19, $t_{221,233,277,391}$ | 1 | 1.289926e+04, $t_{455}$ | 0.00055, $t_{230}$ |
| River-1 | 2.9668, $t_4$ | 3, $t_4$ | 1, $t_4$ | 5, $t_4$ | 0.20, $t_4$ |

*No interactions at 2-simplex level.
**No higher order interaction found.

**Table 20: Connectivity of networks at different simplex levels**

|  | **0-simplex level** | **1-simplex level** | **2-simplex level** |
|---|---|---|---|
| Forest-2 | Disconnected | Connected | Disconnected |
| Forest-4 | Connected | Connected | Disconnected |
| Lake-1 | Connected | Connected | Disconnected |
| Lake-3 | Disconnected | Disconnected | Disconnected |
| Caribbean Food Web network | Connected | Connected | Connected |
| Estuary-2 | Connected | Connected | Connected |
| Rocky shore-1 | Connected | Connected | NA |
| Rocky shore-4 | Connected | Connected | Connected |
| Marine-1 | Connected | Connected | Connected |
| River-1 | Connected | Connected | Connected |

## 5 Discussion and conclusion

Real-world systems may or may not have solely dyadic or paired relations among its vertices; nonetheless, the system's relations may be polyadic in nature, involving more than two vertices at the same time. The clique complex of all the networks has been constructed for better understanding of higher order interactions of carbon mediation networks. For the construction of corresponding clique complexes, the adjacency rules have been redefined by which the interaction of 1-simplices and 2-simplices at 1-simplex and 2-simplex level respectively can be understood. In these clique complexes a 1-simplex represents the relation of two ecological units which are involved in the carbon mediation and for better understanding the significance of simplicial adjacency rules the simplicial complex of Caribbean food-web network can be seen (Fig.7), in which two 1-simplices are adjacent if they have a common 0-simplex and not part of a 2-simplex in the network of Caribbean Food web. It means two 1-simplices are adjacent if they have a common source for carbon mediation and not involved in the formation of 2-simplex of carbon mediation. This result shows the significance of the adjacency at the 1-simplex level. Similarly, adjacency at 2-simplex level can also be understood.

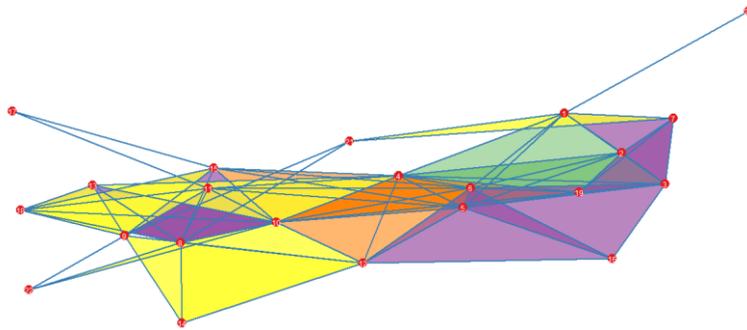

Fig.12: Illustration of the simplicial complex of the Lake 3 network from aquatic environment. Red dots represents the 0-simplices, blue lines are 1-simplices, yellow triangles are the 2-simplices, purple quadrilaterals are the 3-simplices, orange pentagons are 4-simplices and green hexagons are 5-simplices.

After that centrality indices of all the clique complexes have been calculated by extending the definition of graph centrality indices to simplicial centrality indices. Simplicial centrality measures like degree centrality, betweenness centrality, eigenvector centrality, closeness centrality and subgraph centrality have been calculated up to 2-simplex level. Centrality measures calculations showed significant difference in the highest values at both levels. By the centralities at 0-simplex level of the Caribbean Food Web network it can be seen that the vertex $v_1$ has the highest score for all centrality indices, which concluded that vertex $v_1$ is the ecological unit which involved in maximum transfer of the carbon with other ecological units. And it can be concluded that $v_1$ is the most important vertex in this network. But after constructing the corresponding clique complex at 1-simplex and 2-simplex level, we identified that the top ranked *k*-simplices are those who involved in higher-order interaction with $v_1$. Hence, it can be concluded that $v_1$ participated in lower order interaction as well as in the higher-order interaction simultaneously. $v_1$ remains a very highly central vertex in both, the 1-skeleton as also the corresponding clique complex.

The low centrality scores of some 0-simplices corresponding to top ranked 1-simplices and 2-simplices shows less involvement in the carbon mediation at the 0-simplex level. By calculating the centralities scores at different levels, it can be concluded that the bottom ranked 0-simplices can also participate in top ranked simplices at different level. Hence some ecological units which are not significantly involved in the pairwise interaction can participate in higher order interactions.

For example, the betweenness centrality of a 0-simplex $\{v_2\}$ in estuary-2 network at 0-simplex level is 0 while the corresponding 1-simplex $\{v_2, v_{25}\}$ at 1-simplex level is 8.1326. The 0-simplex $\{v_{25}\}$ which is the most valuable 0-simplex of this network, participated with the least valuable 0-simplex $\{v_2\}$ and transformed it into recognisable 1-simplex at 1-simplex level. Similarly, in the CFW network bottom ranked 0-simplices $\{v_{29}\}$ and $\{v_{32}\}$ participated with $v_1$ and $v_4$ to form valuable 2-simplices $\{v_{29}, v_1, v_4\}$ and $\{v_{32}, v_1, v_4\}$ respectively. In carbon mediation networks, a vertex represents an ecological unit (biotic as well as a biotic), and it may be concluded that at 0-simplex level some ecological units are valuable because they are dominantly participating in carbon mediation. But at the 1-simplex level same ecological unit become less or more valuable, in comparison of the group of ecological units which are more dominantly participating in carbon mediation.

The different graph theoretic centralities and simplicial centralities for all the networks have been plotted, and it may be concluded that the ranking of 0,1 and 2-simplices according to a given centrality measure can differ significantly for every simplicial complex which is independent from the region of origin of corresponding networks. On the other hand, the simplicial centrality plot for every network shows that top ranked 0-simplices appears in top ranked 1-simplices and 2-simplices and it may be concluded that 0-simplices with high centrality scores are likely to participate in the higher order interactions. However, some bottom ranked 0-simplices can also participate in higher order interactions.

In the Lake-1 network ranking of first four 0-simplices with respect to degree centrality is $v_{19}, v_2, v_{12}, v_{10}$ while according to closeness centrality the ranking is $v_{19}, v_2, v_{12}, v_{13}$ and nearly same pattern for betweenness centrality as well. However, ranking of 1-simplices and 2-simplices in Lake-1 network is also showing the same variation with respect to same centrality indexes. Simplicial centrality plots for River-1 network from the aquatic environment and Caribbean Food Web network from terrestrial environment of 0-simplices, 1-simplices and 2-simplices with respect to the degree, closeness and betweenness centrality is given in the Plot 1 and Plot 2 respectively.

Forest-2 and Lake-3 are disconnected networks and for the construction of corresponding clique complex the isolated 0-simplex has been deleted. The connectivity of corresponding clique complexes shows difference at 1-simplex level and the with the help of Theorem 3 it can be concluded that clique complex of Forest-2 network is connected at 1-simplex while clique complex of Forest-2 network is disconnected at 1-simplex level. This result suggests us understanding of connectivity of a network at one level does not imply the understanding at another level. Clique complex of any connected or disconnected network can be connected or disconnected at another level, connectivity of all networks at different simplex levels are given in Table 20.

As demonstrated in this work, Definition 12 holds a promise to facilitate a deeper insight into and better the understanding of the structural intricacies of the ecosystem that such a network represents. We have presented arguments for modelling the selected ecosystems in consonance with the definition, with an aim to consequently enhance the understanding of higher-order interactions of carbon mediation networks via the comparison of centrality indices at different simplex levels. Thus, the important higher-order interactions of a given network originating either from aquatic or terrestrial environment, with the ecological units which might be responsible for these non-pairwise interaction can be identified. Construction of clique complex for each of these networks gives us better insight into the carbon mediation networks, and a more detailed evaluation of the structural intricacies in terms of the structural indices have been obtained. Our modelling addresses the first of the two challenges mentioned in the Introduction. Though we do not discuss the second challenge in this work, however, our computations and observations on the higher-order structures of these networks may prove helpful for researchers to explore and address the point in detail.

# References

[1] Levin SA. Ecosystems and the biosphere as complex adaptive systems. Ecosystems. 1(5):431-6; 1998.

[2] Levin S, Xepapadeas T, Crépin AS, Norberg J, De Zeeuw A, Folke C, Hughes T, Arrow K, Barrett S, Daily G, Ehrlich P. Social-ecological systems as complex adaptive systems: modeling and policy implications. Environment and Development Economics. 18(2):111-32; 2013.

[3] Levin SA, Clark W. Toward a science of sustainability. Center for International Development Working Papers; 2010.

[4] Bodin Ö, Norberg J. A network approach for analysing spatially structured populations in fragmented landscape. Landscape Ecology. 22(1):31-44; 2007.

[5] Bell G. The distribution of abundance in neutral communities. The American Naturalist. 155(5):606-17; 2000.

[6] Jordán F, Scheuring I. Network ecology: topological constraints on ecosystem dynamics. Physics of Life Reviews. 1(3):139-72; 2004.

[7] Abrams PA, Matsuda H. Positive indirect effects between prey species that share predators. Ecology 77:610-616; 1996.

[8] Kéfi S, Berlow EL, Wieters EA, Joppa LN, Wood SA, Brose U, Navarrete SA. Network structure beyond food webs: mapping non-trophic and trophic interactions on Chilean rocky shores. Ecology. 96(1):291-303; 2015.

[9] Bascompte J. Networks in ecology. Basic and Applied Ecology. 8(6):485-90; 2007.

[10] Proulx SR, Promislow DE, Phillips PC. Network thinking in ecology and evolution. Trends in ecology & evolution. 20(6):345-53; 2005.

[11] Ings TC, Montoya JM, Bascompte J, Blüthgen N, Brown L, Dormann CF, Edwards F, Figueroa D, Jacob U, Jones JI, Lauridsen RB. Ecological networks–beyond food webs. Journal of animal ecology. 78(1):253-69; 2009.

[12] Jordán F, Scheuring I. Network ecology: topological constraints on ecosystem dynamics. Physics of Life Reviews. 1(3):139-72; 2004.

[13] Higashi M, editor. Theoretical studies of ecosystems: the network perspective. Cambridge University Press; 1991.

[14] Montoya JM, Rodríguez MA, Hawkins BA. Food web complexity and higher-level ecosystem services. Ecology letters. 6(7):587-93; 2003.

[15] Paine RT. A note on trophic complexity and community stability. The American Naturalist. 103(929):91-3; 1969.

[16] Rakshit N, Banerjee A, Mukherjee J, Chakrabarty M, Borrett SR, Ray S. Comparative study of food webs from two different time periods of Hooghly Matla estuarine system, India through network analysis. Ecological Modelling. 356:25-37; 2017.

[17]   Buchary EA, Alder J, Nurhakim S, Wagey T. The use of ecosystem-based modelling to investigate multi-species management strategies for capture fisheries in the Bali Strait, Indonesia. Fisheries Centre Research Reports. 10(2):24; 2002.

[18]   Christensen V, Pauly D, editors. Trophic models of aquatic ecosystems. World Fish; 1993.

[19]   Neira S, Arancibia H, Cubillos L. Comparative analysis of trophic structure of commercial fishery species off Central Chile in 1992 and 1998. Ecological Modelling. 172(2-4):233-48; 2004.

[20]   Schmitz OJ, Leroux SJ. Food webs and ecosystems: linking species interactions to the carbon cycle. Annual Review of Ecology, Evolution, and Systematics. 51:271-95; 2020.

[21]   Roy U, Sarwardi S, Majee NC, Ray S. Effect of salinity and fish predation on zooplankton dynamics in Hooghly–Matla estuarine system, India. Ecological informatics. 35:19-28; 2016.

[22]   Upadhyay S, Roy A, Ramprakash M, Idiculla J, Kumar AS, Bhattacharya S. A network theoretic study of ecological connectivity in Western Himalayas. Ecological Modelling. 359:246-57; 2017.

[23]   Bascompte J. Structure and dynamics of ecological networks. Science. 329(5993):765-6; 2010.

[24]   Pascual M, Dunne JA, Dunne JA, editors. Ecological networks: linking structure to dynamics in food webs. Oxford University Press; 2006.

[25]   Bascompte J. Structure and dynamics of ecological networks. Science. 329(5993):765-6; 2010.

[26]   Song C, Saavedra S. Telling ecological networks apart by their structure: An environment-dependent approach. PLoS computational biology. 16(4):e1007787; 2020.

[27]   Estrada E. The structure of complex networks: theory and applications. Oxford University Press; 2012.

[28]   Bhattacharya S. Higher-order social-ecological network as a simplicial complex. arXiv preprint arXiv:2112.15318; 2021.

[29]   Pascual M, Dunne JA, Dunne JA, editors. Ecological networks: linking structure to dynamics in food webs. Oxford University Press; 2006.

[30]   Dunne JA, Williams RJ, Martinez ND. Food-web structure and network theory: the role of connectance and size. Proceedings of the National Academy of Sciences. 99(20):12917-22; 2002.

[31]   Battiston F, Cencetti G, Iacopini I, Latora V, Lucas M, Patania A, Young JG, Petri G. Networks beyond pairwise interactions: structure and dynamics. Physics Reports. 874:1-92; 2020.

[32]   Mayfield MM, Stouffer DB. Higher-order interactions capture unexplained complexity in diverse communities. Nature ecology & evolution. 1(3):1-7; 2017.

[33] Lang JM, Benbow ME. Species interactions and competition. Nature education knowledge. 4(4):8; 2013.

[34] Huie JM, Prates I, Bell RC, de Queiroz K. Convergent patterns of adaptive radiation between island and mainland Anolis lizards. Biological Journal of the Linnean Society. 134(1):85-110; 2021.

[35] Upadhyay S, Bhattacharya S. A spectral graph theoretic study of predator-prey networks. arXiv preprint arXiv:1901.02883; 2019.

[36] Benson AR, Abebe R, Schaub MT, Jadbabaie A, Kleinberg J. Simplicial closure and higher-order link prediction. Proceedings of the National Academy of Sciences. (48):E11221-30; 2018.

[37] Aktas ME, Nguyen T, Jawaid S, Riza R, Akbas E. Identifying critical higher-order interactions in complex networks. Scientific reports. 11(1):1-1; 2021.

[38] Williams RJ, Berlow EL, Dunne JA, Barabási AL, Martinez ND. Two degrees of separation in complex food webs. Proceedings of the National Academy of Sciences. 99(20):12913-6; 2002.

[39] Pascual M, Dunne JA, Dunne JA, editors. Ecological networks: linking structure to dynamics in food webs. Oxford University Press; 2006.

[40] https://www.globalwebdb.com/

[41] Weinreich DM, Lan Y, Wylie CS, Heckendorn RB. Should evolutionary geneticists worry about higher-order epistasis? Current opinion in genetics & development. 2013 Dec 1;23(6):700-7; 2013.

[42] Mayfield MM, Stouffer DB. Higher-order interactions capture unexplained complexity in diverse communities. Nature ecology & evolution. 1(3):1-7; 2017.

[43] Horak D, Maletić S, Rajković M. Persistent homology of complex networks. Journal of Statistical Mechanics: Theory and Experiment. 2009(03):P03034; 2009.

[44] Jonsson J. Simplicial complexes of graphs. Springer Science & Business Media; 2007.

[45] Van Dalen D, van Dalen D. Logic and structure. Berlin: Springer; 1994.

[46] Shoenfield JR. Mathematical Logic. Addison-Wesley, Reading; 1967.

[47] Srivastava SM, Gödel K. A course on mathematical logic. Springer; 2008.

[48] Battiston F, Cencetti G, Iacopini I, Latora V, Lucas M, Patania A, Young JG, Petri G. Networks beyond pairwise interactions: structure and dynamics. Physics Reports. 874:1-92; 2020.

[49] Torres L, Blevins AS, Bassett D, Eliassi-Rad T. The why, how, and when of representations for complex systems. SIAM Review.63(3):435-85; 2021.

[50] Maunder CR. Algebraic topology. Courier Corporation; 1996.

[51] Hatcher A. Algebraic Topology. Cambridge University Press; 2002.

[52] Deo S. Algebraic Topology: A Primer. Springer; 2018.

[53] Estrada E, Ross GJ. Centralities in simplicial complexes. Applications to protein interaction networks. Journal of theoretical biology. 438:46-60; 2018.

[54] Serrano DH, Gómez DS. Centrality measures in simplicial complexes: Applications of topological data analysis to network science. Applied Mathematics and Computation. 382:125331; 2020.

[55] Brandes U. Network analysis: Methodological foundations. Springer Science & Business Media; 2005.

[56] Scotti M, Jordán F. Relationships between centrality indices and trophic levels in food webs. Community Ecology. 11(1):59-67; 2010.

[57] Pilosof S, Porter MA, Pascual M, Kéfi S. The multilayer nature of ecological networks. Nat Ecol Evol. 1(4):101; 2017.

[58] Grigor'yan A, Muranov YV, Yau ST. Graphs associated with simplicial complexes. Homology, Homotopy and Applications. 16(1):295-311; 2014.

[59] Eaton KA. The life history and production of Chaoborus punctipennis (Diptera: Chaoboridae) in lake Norman, North Carolina, USA. Hydrobiologia. 106(3):247-52; 1983.

[60] Yan ND, Keller W, MacIsaac HJ, McEachern LJ. Regulation of zooplankton community structure of an acidified lake by Chaoborus. Ecological Applications. 1(1):52-65; 1991.

[61] Pascual M, Dunne JA, Dunne JA, editors. Ecological networks: linking structure to dynamics in food webs. Oxford University Press; 2006.

[62] Pedersen AB, Fenton A. Emphasizing the ecology in parasite community ecology. Trends in ecology & evolution. 22(3):133-9; 2007.

[63] Burkle LA, Marlin JC, Knight TM. Plant-pollinator interactions over 120 years: loss of species, co-occurrence, and function. Science. 339(6127):1611-5; 2013.

[64] Estrada E, Rodriguez-Velazquez JA. Subgraph centrality in complex networks. Physical Review E. 71(5):056103; 2005.

[65] Stevanović D. Comment on "Subgraph centrality in complex networks". Physical Review E. 88(2):026801; 2013.

[66] Christian RR, Luczkovich JJ. Organizing and understanding a winter's seagrass foodweb network through effective trophic levels. Ecological modelling. 117(1):99-124; 1999.

[67] Stewart TJ, Sprules WG. Carbon-based balanced trophic structure and flows in the offshore Lake Ontario food web before (1987–1991) and after (2001–2005) invasion-induced ecosystem change. Ecological Modelling. 222(3):692-708; 2011.

[68] Yodzis P. Local trophodynamics and the interaction of marine mammals and fisheries in the Benguela ecosystem. Journal of Animal Ecology. 67(4):635-58; 1998.

[69] Hartley PH. Food and feeding relationships in a community of fresh-water fishes. The Journal of Animal Ecology. 1:1-4; 1948.

[70] Torres MÁ, Coll M, Heymans JJ, Christensen V, Sobrino I. Food-web structure of and fishing impacts on the Gulf of Cadiz ecosystem (South-western Spain). Ecological modelling. 265:26-44; 2013.

[71] Boit A, Martinez ND, Williams RJ, Gaedke U. Mechanistic theory and modelling of complex food-web dynamics in Lake Constance. Ecology letters. 15(6):594-602; 2012.

[72] Angelini R, Agostinho AA. Food web model of the Upper Paraná River Floodplain: description and aggregation effects. Ecological Modelling. 181(2-3):109-21; 2005.

[73] Rasmussen DI. Biotic communities of Kaibab plateau, Arizona. Ecological Monographs. 11(3):230-75; 1941.

[74] Twomey AC. The bird population of an elm-maple forest with special reference to aspection, territorialism, and coactions. Ecological Monographs. 15(2):173-205;1945.

[75] Goldwasser L, Roughgarden J. Construction and Analysis of a Large Caribbean Food Web: Ecological Archives E074-001. Ecology. 74(4):1216-33; 1993.